\catcode`\@=11

\font\tenmsa=msam10 \font\sevenmsa=msam7 \font\fivemsa=msam5
\font\tenmsb=msbm10
\font\sevenmsb=msbm7 \font\fivemsb=msbm5 \newfam\msafam \newfam\msbfam
\textfont\msafam=\tenmsa \scriptfont\msafam=\sevenmsa
\scriptscriptfont\msafam=\fivemsa \textfont\msbfam=\tenmsb
\scriptfont\msbfam=\sevenmsb \scriptscriptfont\msbfam=\fivemsb

\def\hexnumber@#1{\ifnum#1<10 \number#1\else \ifnum#1=10 A\else\ifnum#1=11
 B\else\ifnum#1=12 C\else \ifnum#1=13 D\else\ifnum#1=14 E\else\ifnum#1=15
 F\fi\fi\fi\fi\fi\fi\fi}

\def\msa@{\hexnumber@\msafam} \def\msb@{\hexnumber@\msbfam}
\mathchardef\boxdot="2\msa@00 \mathchardef\boxplus="2\msa@01
\mathchardef\boxtimes="2\msa@02 \mathchardef\square="0\msa@03
\mathchardef\blacksquare="0\msa@04 \mathchardef\centerdot="2\msa@05
\mathchardef\lozenge="0\msa@06 \mathchardef\blacklozenge="0\msa@07
\mathchardef\circlearrowright="3\msa@08 \mathchardef\circlearrowleft="3\msa@09
\mathchardef\rightleftharpoons="3\msa@0A
\mathchardef\leftrightharpoons="3\msa@0B \mathchardef\boxminus="2\msa@0C
\mathchardef\Vdash="3\msa@0D \mathchardef\Vvdash="3\msa@0E
\mathchardef\vDash="3\msa@0F \mathchardef\twoheadrightarrow="3\msa@10
\mathchardef\twoheadleftarrow="3\msa@11 \mathchardef\leftleftarrows="3\msa@12
\mathchardef\rightrightarrows="3\msa@13 \mathchardef\upuparrows="3\msa@14
\mathchardef\downdownarrows="3\msa@15 \mathchardef\upharpoonright="3\msa@16
 \mathchardef\downharpoonright="3\msa@17
\mathchardef\upharpoonleft="3\msa@18 \mathchardef\downharpoonleft="3\msa@19
\mathchardef\rightarrowtail="3\msa@1A \mathchardef\leftarrowtail="3\msa@1B
\mathchardef\leftrightarrows="3\msa@1C \mathchardef\rightleftarrows="3\msa@1D
\mathchardef\Lsh="3\msa@1E \mathchardef\Rsh="3\msa@1F
\mathchardef\rightsquigarrow="3\msa@20
\mathchardef\leftrightsquigarrow="3\msa@21 \mathchardef\looparrowleft="3\msa@22
\mathchardef\looparrowright="3\msa@23 \mathchardef\circeq="3\msa@24
\mathchardef\succsim="3\msa@25 \mathchardef\gtrsim="3\msa@26
\mathchardef\gtrapprox="3\msa@27 \mathchardef\multimap="3\msa@28
\mathchardef\therefore="3\msa@29 \mathchardef\because="3\msa@2A
\mathchardef\doteqdot="3\msa@2B 
\mathchardef\traceiangleq="3\msa@2C \mathchardef\precsim="3\msa@2D
\mathchardef\lesssim="3\msa@2E \mathchardef\lessapprox="3\msa@2F
\mathchardef\eqslantless="3\msa@30 \mathchardef\eqslantgtr="3\msa@31
\mathchardef\curlyeqprec="3\msa@32 \mathchardef\curlyeqsucc="3\msa@33
\mathchardef\preccurlyeq="3\msa@34 \mathchardef\leqq="3\msa@35
\mathchardef\leqslant="3\msa@36 \mathchardef\lessgtr="3\msa@37
\mathchardef\backprime="0\msa@38 \mathchardef\risingdotseq="3\msa@3A
\mathchardef\fallingdotseq="3\msa@3B \mathchardef\succcurlyeq="3\msa@3C
\mathchardef\geqq="3\msa@3D \mathchardef\geqslant="3\msa@3E
\mathchardef\gtrless="3\msa@3F \mathchardef\sqsubset="3\msa@40
\mathchardef\sqsupset="3\msa@41
\mathchardef\trianglelefteq="3\msa@45 \mathchardef\bigstar="0\msa@46
\mathchardef\between="3\msa@47 \mathchardef\blacktriangledown="0\msa@48
\mathchardef\blacktriangleright="3\msa@49
\mathchardef\blacktriangleleft="3\msa@4A
\mathchardef\blacktriangle="0\msa@4E \mathchardef\triangledown="0\msa@4F
\mathchardef\eqcirc="3\msa@50 \mathchardef\lesseqgtr="3\msa@51
\mathchardef\gtreqless="3\msa@52 \mathchardef\lesseqqgtr="3\msa@53
\mathchardef\gtreqqless="3\msa@54 \mathchardef\Rrightarrow="3\msa@56
\mathchardef\Lleftarrow="3\msa@57 \mathchardef\veebar="2\msa@59
\mathchardef\barwedge="2\msa@5A \mathchardef\doublebarwedge="2\msa@5B
\mathchardef\angle="0\msa@5C \mathchardef\measuredangle="0\msa@5D
\mathchardef\sphericalangle="0\msa@5E \mathchardef\varpropto="3\msa@5F
\mathchardef\smallsmile="3\msa@60 \mathchardef\smallfrown="3\msa@61
\mathchardef\Subset="3\msa@62 \mathchardef\Supset="3\msa@63
\mathchardef\Cup="2\msa@64  \mathchardef\Cap="2\msa@65
 \mathchardef\curlywedge="2\msa@66
\mathchardef\curlyvee="2\msa@67 \mathchardef\leftthreetimes="2\msa@68
\mathchardef\rightthreetimes="2\msa@69 \mathchardef\subseteqq="3\msa@6A
\mathchardef\supseteqq="3\msa@6B \mathchardef\bumpeq="3\msa@6C
\mathchardef\Bumpeq="3\msa@6D \mathchardef\lll="3\msa@6E 
\mathchardef\ggg="3\msa@6F  \mathchardef\circledS="0\msa@73
\mathchardef\pitchfork="3\msa@74 \mathchardef\dotplus="2\msa@75
\mathchardef\backsim="3\msa@76 \mathchardef\backsimeq="3\msa@77
\mathchardef\complement="0\msa@7B \mathchardef\intercal="2\msa@7C
\mathchardef\circledcirc="2\msa@7D \mathchardef\circledast="2\msa@7E
\mathchardef\circleddash="2\msa@7F \def\ulcorner{\delimiter"4\msa@70\msa@70 }
\def\urcorner{\delimiter"5\msa@71\msa@71 }
\def\llcorner{\delimiter"4\msa@78\msa@78 }
\def\lrcorner{\delimiter"5\msa@79\msa@79 } \def\yen{\mathhexbox\msa@55 }
\def\checkmark{\mathhexbox\msa@58 } \def\circledR{\mathhexbox\msa@72 }
\def\maltese{\mathhexbox\msa@7A } \mathchardef\lvertneqq="3\msb@00
\mathchardef\gvertneqq="3\msb@01 \mathchardef\nleq="3\msb@02
\mathchardef\ngeq="3\msb@03 \mathchardef\nless="3\msb@04
\mathchardef\ngtr="3\msb@05 \mathchardef\nprec="3\msb@06
\mathchardef\nsucc="3\msb@07 \mathchardef\lneqq="3\msb@08
\mathchardef\gneqq="3\msb@09 \mathchardef\nleqslant="3\msb@0A
\mathchardef\ngeqslant="3\msb@0B \mathchardef\lneq="3\msb@0C
\mathchardef\gneq="3\msb@0D \mathchardef\npreceq="3\msb@0E
\mathchardef\nsucceq="3\msb@0F \mathchardef\precnsim="3\msb@10
\mathchardef\succnsim="3\msb@11 \mathchardef\lnsim="3\msb@12
\mathchardef\gnsim="3\msb@13 \mathchardef\nleqq="3\msb@14
\mathchardef\ngeqq="3\msb@15 \mathchardef\precneqq="3\msb@16
\mathchardef\succneqq="3\msb@17 \mathchardef\precnapprox="3\msb@18
\mathchardef\succnapprox="3\msb@19 \mathchardef\lnapprox="3\msb@1A
\mathchardef\gnapprox="3\msb@1B \mathchardef\nsim="3\msb@1C
\mathchardef\napprox="3\msb@1D
\mathchardef\nsupseteqq="3\msb@23 \mathchardef\subsetneqq="3\msb@24
\mathchardef\supsetneqq="3\msb@25
\mathchardef\supsetneq="3\msb@29 \mathchardef\nsubseteq="3\msb@2A
\mathchardef\nsupseteq="3\msb@2B \mathchardef\nparallel="3\msb@2C
\mathchardef\nmid="3\msb@2D \mathchardef\nshortmid="3\msb@2E
\mathchardef\nshortparallel="3\msb@2F \mathchardef\nvdash="3\msb@30
\mathchardef\nVdash="3\msb@31 \mathchardef\nvDash="3\msb@32
\mathchardef\nVDash="3\msb@33 \mathchardef\ntrianglerighteq="3\msb@34
\mathchardef\ntrianglelefteq="3\msb@35 \mathchardef\ntriangleleft="3\msb@36
\mathchardef\ntriangleright="3\msb@37 \mathchardef\nleftarrow="3\msb@38
\mathchardef\nrightarrow="3\msb@39 \mathchardef\nLeftarrow="3\msb@3A
\mathchardef\nRightarrow="3\msb@3B \mathchardef\nLeftrightarrow="3\msb@3C
\mathchardef\nleftrightarrow="3\msb@3D \mathchardef\divideontimes="2\msb@3E
\mathchardef\varnothing="0\msb@3F \mathchardef\nexists="0\msb@40
\mathchardef\mho="0\msb@66 \mathchardef\thorn="0\msb@67
\mathchardef\beth="0\msb@69 \mathchardef\gimel="0\msb@6A
\mathchardef\daleth="0\msb@6B \mathchardef\lessdot="3\msb@6C
\mathchardef\gtrdot="3\msb@6D \mathchardef\ltimes="2\msb@6E
\mathchardef\rtimes="2\msb@6F \mathchardef\shortmid="3\msb@70
\mathchardef\shortparallel="3\msb@71 \mathchardef\smallsetminus="2\msb@72
\mathchardef\thicksim="3\msb@73 \mathchardef\thickapprox="3\msb@74
\mathchardef\approxeq="3\msb@75 \mathchardef\succapprox="3\msb@76
\mathchardef\precapprox="3\msb@77 \mathchardef\curvearrowleft="3\msb@78
\mathchardef\curvearrowright="3\msb@79 \mathchardef\digamma="0\msb@7A
\mathchardef\varkappa="0\msb@7B \mathchardef\hslash="0\msb@7D
\mathchardef\hbar="0\msb@7E \mathchardef\backepsilon="3\msb@7F
\def\Bbb{\ifmmode\let\next\Bbb@\else
\def\next{\errmessage{Use \string\Bbb\space only in math mode}}\fi\next}
\def\Bbb@#1{{\Bbb@@{#1}}} \def\Bbb@@#1{\fam\msbfam#1}

\catcode`\@=\active

\def\Z{{\Bbb Z}}

\def\lform{\hbox{$\sqcup$}\llap{\hbox{$\sqcap$}}}

\def\eps{{\epsilon}}

\def\<{\langle}
\def\>{\rangle}

\def\lcross{{>\!\!\!\triangleleft}}

\def\cobicross{{\triangleright\!\!\!\blacktriangleleft}}
\def\bicross{{\blacktriangleright\!\!\!\triangleleft}}

\def\tens{\mathop{\otimes}}
\def\la{{\triangleright}}\def\ra{{\triangleleft}}

\def\id{{\rm id}}

\def\o{{}_{(1)}}\def\t{{}_{(2)}}\def\th{{}_{(3)}}

\def\proof{\goodbreak\noindent{\bf Proof\quad}}
\def\endproof{{\ $\lform$}\bigskip }

\def\note#1{}

\def\eqn#1#2{\begin{equation}#2\label{#1}\end{equation}}

\def\text#1{{\rm #1}}

\documentstyle[11pt]{article}
\newtheorem{lemma}{Lemma}[section]
\newtheorem{propos}[lemma]{Proposition}
\newtheorem{example}[lemma]{Example}

\newtheorem{corol}[lemma]{Corollary}

\newtheorem{remark}[lemma]{Remark}

\begin{document}
\baselineskip 22pt

\newpage

{\ }\hskip 4.7in DAMTP/94-11 
\vspace{.5in}

\begin{center} {\Large FINITE GROUP FACTORIZATIONS AND BRAIDING}
\\ \baselineskip 13pt{\ }
{\ }\\ E. J. Beggs \& J. D. Gould\\
{\ } \\ Department of Mathematics\\
University College of Swansea, Swansea\\ Wales SA2 8PP\\
{\ } \\
\& \\
{\ }\\
S. Majid\footnote{Royal Society Univ. Research Fellow and Fellow
 of Pembroke
College, Cambridge.
On leave at Dept Mathematics, Harvard University
during 1995+1996.}\\ {\ }\\
Department of Applied Mathematics \& Theoretical Physics\\ University of
Cambridge, Cambridge CB3 9EW, U.K.
\end{center}
February 1994 -- revised March 1995
\begin{center}
\end{center}
\vspace{10pt}
\begin{quote}\baselineskip 13pt
\noindent{\bf ABSTRACT} We compute the quantum double, braiding and other
canonical Hopf algebra
constructions for the bicrossproduct Hopf algebra $H$ associated to the
factorization of a finite group into
two subgroups. The representations of the quantum double are described
by a notion of
bicrossed bimodules, generalising the cross modules of Whitehead. We also
show
that self-duality structures for the bicrossproduct Hopf algebras are in
one-one correspondence with factor-reversing group isomorphisms. The
example
$Z_6Z_6$ is given in detail.
 We show further that the quantum double $D(H)$ is the twisting of
$D(X)$ by a non-trivial quantum cocycle, where $X$ is the associated
double cross product group.
\end{quote}
\baselineskip 22pt

\section{INTRODUCTION}

It is known that to every factorisation $X=GM$ of a group into two subgroups
$G,M$, is associated a generally non-commutative and non-cocommutative Hopf
algebra $H=kM\cobicross k(G)$. Here $k(G)$ is the
Hopf algebra of functions on $G$ and $kM$ is the group Hopf algebra of $M$.
Further details will be recalled in the Preliminaries. These bicrossproduct
Hopf algebras arose in an algebraic approach to quantum-gravity in
\cite{Ma:phy} as well as having been noted in connection
with extension theory in \cite{Tak:mat}. Group factorisations are very common
in mathematics and this bicrossproduct  construction remains one
 of the
primary sources of true non-commutative and non-cocommutative Hopf
algebras. One
of their novel features, which was the motivation in \cite{Ma:phy} is  the
self-dual nature
of the construction. The dual Hopf algebra is of the same general form, namely
$k(M)\bicross kG$.

Since this work, there have been developed a number of more modern Hopf
algebra
constructions related to
knot and three manifold invariants. Central to these is the Quantum Double
construction of V.G. Drinfeld\cite{Dri} which associates to a general Hopf
algebra $H$ a quasitriangular one $D(H)$. This
quantum double Hopf algebra induces on its category of representations a
braiding. When $H=kG$ the group algebra of a group, the braiding is that
of the category of crossed $G$-modules as studied by
Whitehead\cite{Whi:com}.
In the present paper we compute the quantum double and braiding for the
significantly more complicated bicrossproduct Hopf algebras associated to the
factorisation $X=GM$. The result is an interesting generalisation of
crossed modules to bicrossed bimodules.

We also study the question of when exactly the
bicrossproducts $kM\cobicross
k(G)$ are self-dual as Hopf algebras (rather than merely of self-dual form),
and find that such self-duality structures are in one-one
correspondence with factor-reversing automorphisms of $X$. In general one
might hope that many properties of the
bicrossproduct Hopf algebras could be related
to the factorizing group $X$, but this is the first result that we know of this
type. We follow this with a relation between the quantum double of the
bicrossproduct and the double $D(X)$ of the group algebra of $X$.

An outline of the paper is as follows. We begin in Section~2 with the
self-duality result. Section~3 computes the quantum double of the
bicrossproduct Hopf algebra $H=kM\cobicross k(G)$.
Section~4 analyses their representation theory
and computes the canonical Schroedinger representation of $D(H)$ on $H$.
The induced braiding is a Yang-Baxter operator $V\tens V\to V\tens
V$ where $V$ is a vector space with basis $X$. In the paper we compute our
various constructions on one of the simplest examples where the group
factorisation is $\Z_6\Z_6$.
Finally, we show in Section~5 that the braided category of
representations of $D(H)$ as computed in earlier sections, is
monoidally equivalent to the braided category of representations of
$D(X)$ (or crossed X-modules), where $X$ is the double cross product
group. This is disappointing from the point of view of obtaining
completely new braided categories but is an important result from the
algebraic point of view. It means that $D(H)$ is isomorphic as a Hopf
algebra to the twisting of the Hopf algebra $D(X)$ by a Hopf algebra
2-cocycle $F$ in the terminology of \cite{Ma:coh}, \cite{Ma:book}.

The general theory of twisting at the level of quasiHopf algebras was
introduced by Drinfeld in \cite{Dri:2}
 but works also at the Hopf algebra level
\cite{Gur:Maj}. It consists of leaving the algebra
unchanged but conjugating the coproduct by $F$. The sense in which $F$
is a cocycle is itself an interesting one and corresponds to a kind of
quantum or non-Abelian cohomology as explained in
\cite{Ma:coh},
\cite{Ma:book}.  Previous examples have generally been
 Lie or deformation-theoretic:
our results in this section provide $D(H)$ as one
of the first discrete examples. We show that $F$ is non-trivial (not a
coboundary).  It is worth noting that the concept of two cocycle and
2-cohomology here is a special case of the dual of the notion occuring
 in another context in \cite{Doi} and elsewhere.

\subsection*{Acknowledgements}
We would like to thank Alan Thomas for some useful comments.

\subsection*{Preliminaries} For
bicrossproducts we use the conventions of
\cite{Ma:book}. By definition, a group $X$ factorises into $X=GM$
if the map
$G\times M\to X$ given by product in $X$ is a bijection. In this
case define
corresponding action $\la:M\times G\to G$ and
reaction $\ra:M\times G\to M$ by $su=(s\la u)(s\ra u)$
for $s\in M$ and $u\in G$. One can see that they
obey for all $s,t\in M$ and $u,v\in G$:
\[ s\ra e=s,\quad (s\ra u)\ra v=s\ra(uv);\qquad e\ra u=e,\quad (st)\ra
u=\left(s\ra(t\la u)\right)(t\ra u)\]
\[ e\la u=u,\quad s\la(t\la u)=(st)\la u;\quad s\la e=e,\quad s\la(uv)=(s\la
u)\left((s\ra u)\la v\right).\]
One says that $(G,M,\la,\ra)$ are a matched pair of groups acting on each
other. We then use this data to define the
bicrossproduct Hopf algebra
$k(M)\bicross kG$ with basis $\delta_s\tens u$
where $s\in M$ and $u\in G$, and for its dual Hopf algebra
$kM\cobicross k(G)$
with dual basis $s\tens\delta_u$. In both cases they are semidirect
products as
algebras and as coalgebras. Explicit formulae are as follows. For
$kM\cobicross
k(G):$
\[ (s\tens \delta_u)(t\tens \delta_v)=\delta_{u,t\la v}(st\tens
\delta_{v}),\quad \Delta (s\tens \delta_u)=\sum_{xy=u}s\tens\delta_x\tens
s\ra x\tens \delta_y\]
\[ 1=\sum_u e\tens \delta_u,\quad \eps (s\tens \delta_u)=\delta_{u,e},\quad
S(s\tens \delta_u)=(s\ra u)^{-1}\tens\delta_{(s\la u)^{-1}},\quad
(s\tens \delta_u)^*=s^{-1}\tens\delta_{s\la u}
\]
Here $\Delta$ denotes the coproduct, $\eps$ the counit and $S$ the antipode
of
the Hopf algebra. Additionally we have given a formula for the star
operation.
 For $k(M)\bicross kG$:
\[ (\delta_s\tens u)(\delta_t\tens v)=\delta_{s\ra u,t}(\delta_s\tens uv),
\quad
\Delta(\delta_s\tens u)=\sum_{ab=s}\delta_a\tens b\la u\tens
\delta_b\tens u\]
\[ 1=\sum_s\delta_s\tens e,\quad \eps(\delta_s\tens u)=\delta_{s,e},\quad
S(\delta_s\tens u)=\delta_{(s\ra u)^{-1}}\tens(s\la u)^{-1},\quad
(\delta_s\tens u)^*=\delta_{s\ra u}\tens u^{-1}
.\]

This is all that we need from the theory of
bicrossproducts. The general theory
for arbitrary Hopf algebras was developed in \cite{Ma:phy} to which we refer
the interested reader, or see \cite{Ma:book}. Finally, we note that the
factorising group $X$ can also be recovered from the matched pair of groups as
a {\em double cross product} or double-semidirect product construction in
which
both factors act on each other simultaneously. This group $G\bowtie M$
is built
on $G\times M$ and the explicit formulae which we shall need are:
\[ (u,s)(v,t)=(u(s\la v),(s\ra v)t),\quad e=(e,e),\quad (u,s)^{-1}=(s^{-1}\la
u^{-1},s^{-1}\ra u^{-1}).\]
This double cross product group is isomorphic to our original $X$ by
identifying $(u,s)$ with $us$ in $X$.
We will also need double cross products at the level of Hopf algebras [1],
which
general theory we recall at the beginining of the relevant sections. For
general Hopf algebra constructions we use the summation notation
$\Delta h=\sum
h\o\tens h\t$ as in \cite{Swe:hop}.

\section{Self-Duality of bicrossproducts}

Let $X=GM$ be a group factorisation. We define a group isom $\theta:X\to X$ to
be {\em factor-reversing} if
$\theta(G)\subset M$ and $\theta(M)\subset G$.

\begin{propos} Factor-reversing isomorphisms of $X=GM$
give rise to Hopf algebra self-duality pairings $\<\ ,\ \>:H\tens H\to
k$ on the Hopf algebra $H=kM\cobicross k(G)$. The corresponding pairing is
\[ \<s\tens \delta_u,t\tens \delta_v\>=\delta_{s,\theta(t\la
v)}\delta_{u,\theta(t\ra v)}.\]
Further, a Hopf algebra isomorphism $:H\to H^*$
which sends basis elements to basis elements must arise from a
factor-reversing isomorphism of $X=GM$.
\end{propos}
\proof We define a linear map $\tilde\theta:H\to H^*$ by
\[ \tilde\theta(s\tens \delta_u)=\delta_{\theta(s\la u)}\tens\theta(s\ra u)\]
and verify that this is a Hopf algebra isomorphism
$kM\cobicross k(G) \to k(M)\bicross
kG$ if and only if $\theta$ is a group isomorphism.

If $\theta$ is a group homomorphism then
$\theta(tv)=\theta(t\la v)\theta(t\ra v
)$, which is also $\theta(t)\theta(v)$. The condition that these two
expressions are the same is that, for all $t$ and $v$,
\[ \theta(t)=\theta(t\la v)\la \theta(t\ra v) \quad,\quad
\theta(v)=\theta(t\la v)\ra \theta(t\ra v)  \ .\]
On the assumption that $\theta$ is a group isomorphism, we now check the
conditions for $\tilde\theta$ to be an algebra isomorphism.
\[ \tilde\theta\big((s\tens \delta_u)(t\tens \delta_v)\big)\ =\
\delta_{u,t\la v}\ \delta_{\theta(st\la v)}\tens \theta(st\ra v)\ .\]
Calculation the other way gives
\begin{eqnarray*}
 \tilde\theta(s\tens \delta_u) \tilde\theta(t\tens \delta_v)& = &
\delta_{\theta(t\la v),\theta(s\la u)\ra \theta(s\ra u)}\
\delta_{\theta(s\la u)}\tens \theta(s\ra u)\theta(t\ra v)   \\
& = &
\delta_{\theta(t\la v),\theta(u)}\
\delta_{\theta(s\la u)}\tens \theta\big((s\ra u)(t\ra v)\big)   \\
& = &
\delta_{u,t\la v}\
\delta_{\theta(s\la (t\la v))}\tens \theta\big((s\ra(t\la v))(t\ra
v)\big)   \\
& = &
\delta_{u,t\la v}\
\delta_{\theta(st\la v)}\tens \theta\big(st\ra v\big)   \ .
\end{eqnarray*}

We now check the condition for $\tilde\theta$ to be a
coalgebra isomorphism, i.e. $\Delta\tilde\theta=(
\tilde\theta\tens \tilde\theta)\Delta$.
\begin{eqnarray*}
\Delta\tilde\theta(s\tens\delta_u)& = & \sum_{ab=\theta(s\la u)}
\delta_a\tens b\la \theta(s\ra u) \tens \delta_b\tens \theta(s\ra u)  \\
 (\tilde\theta\tens \tilde\theta)\Delta (s\tens\delta_u)& = &
\sum_{xy=u}
\delta_{\theta(s\la x)}\tens  \theta(s\ra x)
\tens \delta_{\theta\big((s\ra x)\la y\big)} \tens \theta(s\ra xy)
\ .\end{eqnarray*}
If we put $a=\theta(s\la x)$ and $b=\theta\big((s\ra x)\la y\big)$, then
$ab=\theta(s\la u)$, and
\begin{eqnarray*}
b\la \theta(s\ra u)& = & \theta\big((s\ra x)\la y\big)\la \theta(s\ra u)  \\
& = & \theta\big((s\la x)^{-1}(s\la u)\big)\la \theta(s\ra u)  \\
& = & \theta(s\la x)^{-1}\la \big( \theta(s\la u)\la \theta(s\ra u) \big) \\
& = & \theta(s\la x)^{-1}\la  \theta(s) \\
& = & \theta(s\la x)^{-1}\la \big( \theta(s\la x)\la \theta(s\ra x) \big) \\
& = & \theta(s\ra x)
\ ,\end{eqnarray*}
so we get a coalgebra map. Next we check the effect of $\tilde\theta$ on the
unit and counit.
\begin{eqnarray*}
\eps_{H^*}\tilde\theta(s\tens\delta_u)& = &
\delta_{\theta(s\la u),e}\ =\ \delta_{u,e}\ =\
\eps_{H}(s\tens\delta_u)\ .  \\
\tilde\theta(1_H)& = & \tilde\theta\big(\sum_{u}e\tens \delta_u\big)\ =\
\sum_u \delta_{\theta(u)}\tens e\ =\ 1_{H^*}\ .
\end{eqnarray*}
To check that the antipode is preserved, we need to note that
\[
  u^{-1}s^{-1}\ =\  (s\ra u)^{-1}(s\la u)^{-1}\ =\
\big((s\ra u)^{-1}\la (s\la u)^{-1}\big)
\big((s\ra u)^{-1}\ra (s\la u)^{-1}\big) \ ,
\]
so that
\[
\theta(u)^{-1}\ =\ \theta\big((s\ra u)^{-1}\la (s\la u)^{-1}\big)
\quad\text{and}\quad
\theta(s)^{-1}\ =\ \theta\big((s\ra u)^{-1}\ra (s\la u)^{-1}\big) \ .
\]
Then we have that
\begin{eqnarray*}
\tilde\theta S (s\tens\delta_u)& = &
\delta_{\theta\big((s\ra u)^{-1}\la (s\la u)^{-1}\big) } \tens
\theta\big((s\ra u)^{-1}\ra (s\la u)^{-1}\big)
\ =\ \delta_{\theta(u)^{-1}}\tens \theta(s)^{-1}  \\
S \tilde\theta  (s\tens\delta_u)& = &
\delta_{\big( \theta(s\la u)\ra \theta(s\ra u)\big)^{-1} }\tens
\big( \theta(s\la u)\la \theta(s\ra u)\big)^{-1}
 \ =\
\delta_{\theta(u)^{-1}}\tens \theta(s)^{-1}
\end{eqnarray*}
Finally to see that $\tilde\theta$ is invertible, put
\[
\tilde\theta^{-1}(\delta_s\tens u)\ =\ \theta^{-1}(s\la u)\tens
\delta_{\theta^{-1}(s\ra u)}\ ,
\]
then
\begin{eqnarray*}
\tilde\theta^{-1}\tilde\theta(s\tens\delta_u)& = &
\tilde\theta^{-1}\big(
\delta_{\theta(s\la u)}\tens \theta(s\ra u)\big)  \\
& = & \theta^{-1}\big( \theta(s\la u)\la \theta(s\ra u) \big)\tens
\delta_{\theta^{-1}\big( \theta(s\la u)\ra \theta(s\ra u) \big)} \\
& = & \theta^{-1}\big( \theta(s) \big)\tens
\delta_{\theta^{-1}\big( \theta( u) \big)}\ =\ s\tens\delta_u \ ,\\
\tilde\theta\tilde\theta^{-1}(\delta_s\tens u)& = &
\delta_{ \theta\big( \theta^{-1}(s\la u)\la \theta^{-1}(s\ra u) \big)} \tens
\theta\big( \theta^{-1}(s\la u)\ra \theta^{-1}(s\ra u) \big)
\ =\ \delta_s\tens u \ .
\end{eqnarray*}
The last line proceeds since we can swap the roles of $\theta$ and
$\theta^{-1}$ in the group identities, since both are isomorphisms.
 This completes the
proof that $\tilde \theta$ is a Hopf algebra isomorphism.

Now we assume that $\tilde \theta$ is a Hopf algebra isomorphism
which sends our basis elements to basis elements for the preferred
desriptions of our two Hopf algebras, and prove that
we can recover  a group isomorphism $\theta$. We start with functions
$m:M\times G\to M$ and $g:M\times G\to G$ such that
\[  \tilde\theta^{-1}(\delta_s\tens u)\ =\ m(s,u)\tens\delta_{g(s,u)}\ .\]
Since the map $m\times g:M\times G\to M\times G$ is a 1-1
correspondence, it is possible to define a map $\theta$ by
\[
(m\times g)^{-1}(s,u)\ =\
\big(\theta(s\la u),\theta(s\ra u)\big)\ .
\]
Now $\tilde\theta$ is invertible, so that
\[
s\tens \delta_u\ =\ \tilde\theta^{-1}\tilde\theta(s\tens \delta_u)
\quad\text{ and }\quad
\delta_s\tens u\ =\ \tilde\theta\tilde\theta^{-1}(\delta_s\tens u)
\ ,\]
giving the relations
\begin{eqnarray*}
s& = &m\big(\theta(s\la u),\theta(s\ra u)\big)\quad (a)\ ,\\
u& = &g\big(\theta(s\la u),\theta(s\ra u)\big)\quad (b)\ ,\\
s& = &\theta\big(m(s, u)\la g(s, u)\big)\quad (c)\ ,\\
u& = &\theta\big(m(s, u)\ra g(s, u)\big)\quad (d)\ .
\end{eqnarray*}
Next we use the fact that $\tilde\theta^{-1}$ is an algebra
homomorphism in the
equation
\[ \tilde\theta^{-1}\big((\delta_s\tens u)(\delta_t\tens v)\big)\ =\
\tilde\theta^{-1}(\delta_s\tens u)\tilde\theta^{-1}(\delta_t\tens v)\]
to get the equation
\[    \delta_{s\ra u,t}\big(m(s,uv)\tens\delta_{g(s,uv)}\big)\ =\
  \delta_{g(s,u),m(t,v)\la g(t,v)}
\big((m(s,u)m(t,v)\tens\delta_{g(t,v)}\big)\ .\]
Thus for all $s,t\in M$ and $u,v\in G$, the following are equivalent:
\begin{eqnarray*}
  t& = &s\ra u\quad (e)\ ,\\
g(s,u)& = &m(t,v)\la g(t,v)\quad (f)\ ,\\
m(s,uv)& = &m(s,u)m(t,v)\quad (g)\ ,\\
g(s,uv)& = &g(t,v)\quad (h)\ .
\end{eqnarray*}
Now define $\psi:X\to X$ by
\[  \psi(s)=g(s,e)\ ,\ \psi(u)=m(e,u)
\ ,\ \psi(su)=g(s,e)m(e,u)=\psi(s)\psi(u)\ .
\]
Note that $\psi$ is well defined since $G\cap M=\{ e \}$. We show that
$\psi\theta=\theta\psi=\id_X$, so that $\theta$ is a factor reversing
isomorphism. From (e) and (f) we have
\[  g(s,u)\ =\ m(s\ra u,v)\la g(s\ra u,v)\ ,\]
and replacing $s$ by $s\ra u^{-1}$ we obtain
\[  g(s\ra u^{-1},u)\ =\ m(s,v)\la g(s,v)\quad (i)\ .\]
The equations (e) and (h) give $g(s,uv)=g(s\ra u,v)$, and
 replacing $s$ by $s\ra u^{-1}$, and setting $v=e$ gives
\[ g(s\ra u^{-1},u)\ =\ g(s,e)\ =\ \psi(s)\quad (j)\ .\]
Then using (c), (i), and (j), we get
 \[
 s\ =\ \theta\big(m(s,u)\la g(s,u)\big)  \ =\ \theta\big(g(s\ra u^{-1},u)\big)\
=\ \theta\psi(s)\ ,
 \]
so that $\theta\psi|_M=\id_M$. Next we note that from (b) and (i),
 \[
 e\ =\ g\big(\theta(e\la e), \theta(e\ra e)\big)\ =\
g\big(\theta(e), \theta(e)\big) \ =\ g\big(e\ra u^{-1},e\big)\ =\ g(e,u)\ .
 \]
Then using (d) with $s=e$, we have
 \[
 u = \theta\big(m(e,u)\ra g(e,u)\big)\ =\ \theta\big(m(e,u)\ra e\big)  \ =\
\theta\big(m(e,u)\big)\ =\ \theta\psi(u)\ ,
 \]
so $\theta\psi|_G=\id_G$. Putting $u=e$ in (a) gives
 \[
 s\ =\  m\big(\theta(s\la e), \theta(s\ra e)\big)\ =\
m\big(\theta(e), \theta(s)\big)  \ =\
 m\big(e, \theta(s)\big) \ =\ \psi\theta(s)\ ,
 \]
so that $\psi\theta|_M=\id_M$. Next putting $s=e$ in (b) gives
 \[
 u\ =\  g\big(\theta(e\la u), \theta(e\ra u)\big)\ =\
g\big(\theta(u), \theta(e)\big) \ =\
\psi\theta(u)\ ,
 \]
so that $\psi\theta|_G=\id_G$. Finally we have from the definition of $\psi$,
\[
\theta\psi(su)\ =\ \theta\big(\psi(s)\psi(u)\big)\ =\
\theta\psi(s)\theta\psi(u)\ =\ su\ ,
\]
and since $\theta(s\la u)\in M$ and $\theta(s\ra u)\in G$,
\[
\psi\theta(su)\ =\ \psi\theta\big((s\la u)(s\ra u)\big)\ =\
(s\la u)(s\ra u)\ =\ su\ .
\]
 \endproof

\begin{example} {\rm
For odd $n\ge 3$, we can express the product of dihedral groups
$D_n\times D_n$ as a
double cross product of cyclic groups $Z_{2n} Z_{2n}$. We give
the $n=3$ case in detail, where we note that $D_3\equiv S_3$.
Consider the group
$X=S_3\times S_3$ as the permutations of 6 objects labelled 1 to 6,
where the
first factor leaves the last 3 objects unchanged, and the second
factor leaves the first 3 objects unchanged. We take $G$ to be the
cyclic group
of order 6 generated by the permutation $1_G=(123)(45)$, and
$M$ to be the cyclic group
of order 6 generated by the permutation $1_M=(12)(456)$.
Our convention is that
permutations act on objects on their right, for example $1_G$
applied to $1$
gives $2$. The intersection of $G$ and $M$ is just the identity
permutation, and
counting elements shows that $GM=MG=S_3\times S_3$. We write
each cyclic group
additively, for example $G=\{0_G,1_G,2_G,3_G,4_G,5_G\}$.
The action of the
element $1_M$ on $G$ is seen to be given by the permutation
$(1_G,5_G)(2_G,4_G)$, and that of  $1_G$ on $M$ is given by the
permutation
$(1_M,5_M)(2_M,4_M)$.
The algebras $H$ and $H^*$ can be given a matrix decomposition,
for example
 \begin{eqnarray*}
H^* &\ \longrightarrow\ & L^1(Z_6)\ \oplus\
L^1(Z_6)\ \oplus\ M_2\big(L^1(Z_3)\big)\ \oplus\
M_2\big(L^1(Z_3)\big)\ ,\\
\sum a_{ij}\delta_{i_M}\otimes j_G&\ \longmapsto\ &
\big(a_{00}\underline 0\ +\ a_{01}\underline 1\ + \dots +\ a_{05}
\underline
5\big)
\ \oplus\
\big(a_{30}\underline 0\ +\ a_{31}\underline 1\ + \dots +\ a_{35}
\underline
5\big)  \\
 & &\oplus\ \left( \begin{array}{cc}
a_{10}\underline 0\ +\ a_{12}\underline 1\ +\ a_{14}\underline 2 &
a_{15}\underline 0\ +\ a_{11}\underline 1\ +\ a_{13}\underline 2 \\
a_{51}\underline 0\ +\ a_{53}\underline 1\ +\ a_{55}\underline 2 &
a_{50}\underline 0\ +\ a_{52}\underline 1\ +\ a_{54}\underline 2
\end{array}\right)
\\ & &\oplus\ \left( \begin{array}{cc}
a_{20}\underline 0\ +\ a_{22}\underline 1\ +\ a_{24}\underline 2 &
a_{25}\underline 0\ +\ a_{21}\underline 1\ +\ a_{23}\underline 2 \\
a_{41}\underline 0\ +\ a_{43}\underline 1\ +\ a_{45}\underline 2 &
a_{40}\underline 0\ +\ a_{42}\underline 1\ +\ a_{44}\underline 2
\end{array}\right)
 \end{eqnarray*}

In this case we have a factor reversing isomorphism on $X$, in fact
several of
them. The map $\theta(x)=\Theta x \Theta^{-1}$ where $\Theta$ is the
permutation
$(1425)(36)$ gives such an isomorphism. Also the map
$\phi(x)=\Phi x \Phi^{-1}$ where $\Phi$ is the permutation
$(14)(25)(36)$ gives such an isomorphism. In fact the dihedral
group $D_4$
generated by $\theta$ (order 4) and $\phi$ (order 2) is all the
isomorphisms
of $X$ which either reverse or preserve the factors.
}
\end{example}

\begin{example}
{\rm
In more generality, take the dihedral group $D_n$ to
be generated by elements $a$ of order $n$, and $b$ of order $2$, with
relation
$bab=a^{-1}$. For $n$ odd, consider the group
$G$ to be the cyclic subgroup of $D_n\times D_n$
of order $2n$ generated by $1_G=(a,b)$, and
$M$ to be the cyclic subgroup
of order $2n$ generated by $1_M=(b,a)$. By counting elements,
we see that everything in $X=D_n\times D_n$ can be written as a product
on $GM$,
and as a product in $MG$. This gives $D_n\times D_n$ as a
double cross product
$\Z_{2n} \Z_{2n}$. The corresponding actions are (writing the cyclic groups
additively)
\[
s\la u\ =\
\left\{ \begin{array}{ll}
 u  & s\ \text{ even}  \\ -u  & s\ \text{ odd}
\end{array} \right.
\quad\text{and}\quad
s\ra u\ =\
\left\{ \begin{array}{ll}
 s  & u\ \text{ even}  \\ -s  & u\ \text{ odd}
\end{array} \right.\ .
\]
There is a factor reversing isomorphism of $X=D_n\times D_n$ given by
$\psi(x,y)=(y,x)$.
}
\end{example}

\begin{example}
{\rm
Let $X$ be a finite group whose order is divisible by only two distinct
primes $p$ and $q$,
i.e.\ $|X|=p^nq^m$. Then Sylow's theorems state that there is at least one
subgroup $G$ with order $p^n$, and  at least one
subgroup $M$ with order $q^m$. Since the order of $G\cap M$ would have to
divide $p^n$ and  $q^m$, we must have $G\cap M=\{ e\}$, giving the
factorisation $X=GM=MG$. However bicrossproducts
formed from factorisations of this sort can never
 be
self-dual, as the orders of $G$ and $M$ are different.
}
\end{example}

\section{Quantum Double of a bicrossproduct}

Let $H$  be a Hopf algebra. The quantum double is a Hopf algebra double cross
product $D(H)=H^{*{\rm op}}\bowtie H$, not to be confused with  the
bicrossproducts
above. It is a Hopf algebra factorising into $H^{*\rm op}$ and $H$ and given
via a double-semidirect product by mutual coadjoint actions of these two
factors on each other. This formulation is from \cite{Ma:phy} although the
double itself as a Hopf algebra originates with V.G. Drinfeld\cite{Dri}.
Explicitly, it is built on $H^*\tens H$ as a linear space with product
\eqn{double}{ (a\tens h)(b\tens g)=\sum \<Sb\o,h\o\> b\t a\tens h\t
g\<b\th,h\th\>}
and tensor product unit, counit, coproduct, and antipode given by the
formulae:
\begin{eqnarray*}
1_{D(H)} & = & 1_{H^*}\tens 1_H  \ , \\
\eps_{D(H)} & = & \eps_{H^*}\tens \eps_H\ ,  \\
\Delta_{D(H)} & = & (\id\tens\tau\tens\id)\circ(\Delta\tens\Delta)\ ,\\
S_{D(H)}(a\tens h) & = & (1\tens Sh)(S^{-1}a\tens 1) \ ,\\
(a\tens h)^* & = & (1\tens h^*)((S^2a)^*\tens 1) \ ,
\end{eqnarray*}
where $\tau(a\tens h)=h\tens a$.

One can compute this directly or compute first the mutual coadjoint actions
$\la,\ra$. In their terms, the Hopf algebra structure takes the double cross
product form
\eqn{dcross}{(a\tens h)(b\tens g)=\sum  (h\o\la b\o)a\tens (h\t\ra b\t)g}
We will use this double cross product form for our calculation. The required
mutual coadjoint actions are defined by
\eqn{coadj}{ h\ra a=\sum h\t \<a,(Sh\o)h\th\>,\quad h\la a=\sum a\t
\<h,(Sa\o)a\th\>,\qquad a\in H^*,\ h\in H.}
In all these formulae, the expresssions are given in terms of the Hopf algebras
$H,H^*$.

\begin{lemma} The coadjoint action of $H=kM\cobicross k(G)$ on
$H^*=k(M)\bicross kG$ and vice-versa are
\[ (t\tens\delta_v)\la (\delta_s\tens u)=\delta_{u,(s\la
u)v}\delta_{t's{t'}^{-1}}\tens t'\la u;\quad t'=t\ra (s\la u)^{-1}\]
\[ (t\tens\delta_v)\ra(\delta_s\tens u)=t'\tens \delta_{(s\la
u)vu^{-1}}\delta_{t\ra v,t(s\ra u)}.\]
\end{lemma}

\proof :  The action of $H$ on $H^*$.
Let $a=\delta_s\tens u\in H$, then firstly we calculate
\[  (\Delta\tens \id)\Delta(\delta_s\tens u)\ =\ \sum_{cdb=s}
\delta_c\tens db\la u\tens \delta_d\tens b\la u\tens\delta_b\tens u\ .\]
{}From this we get that $a_{(1)}=\delta_c\tens db\la u$,
$a_{(2)}=\delta_d\tens b\la u$, and $a_{(3)}=\delta_b\tens u$,  which with
$h=t\tens\delta_v$ we substitute into the second equation of (3) to get
\begin{eqnarray*}
(t\tens\delta_v)\la(\delta_s\tens u)& = &\sum_{cdb=s} (\delta_d\tens b\la u)
\big< t\tens\delta_v,S(\delta_c\tens db\la u)\delta_b\tens u\big>  \\
& = &\sum_{cdb=s} (\delta_d\tens b\la u)
\big< t\tens\delta_v,\big(\delta_{(c\ra(db\la u))^{-1}}\tens (c\la(db\la
u))^{-1}\big) (\delta_b\tens u)\big>
\\
& = &\sum_{cdb=s} (\delta_d\tens b\la u)
\big< t\tens\delta_v,\big(\delta_{c^{-1}\ra(s\la u)}\tens (s\la u)^{-1}\big)
(\delta_b\tens u)\big>  \\
& = &\sum_{cdb=s} (\delta_d\tens b\la u)
\big< t\tens\delta_v,  \delta_{ (c^{-1}\ra(s\la u))\ra (s\la u)^{-1} , b }
(\delta_{c^{-1}\ra(s\la u)}\tens (s\la u)^{-1}u)\big>  \\
& = &\sum_{cdb=s} (\delta_d\tens b\la u)
\big< t\tens\delta_v,  \delta_{ c^{-1} , b }
(\delta_{c^{-1}\ra(s\la u)}\tens (s\la u)^{-1}u)\big>  \\
& = &\sum_{cdb=s} (\delta_d\tens b\la u)
\delta_{ c^{-1} , b } \delta_{ t , c^{-1}\ra(s\la u) } \delta_{ v ,(s\la
u)^{-1}u } \\
& = &\delta_{ v ,(s\la
u)^{-1}u } \sum_{cdc^{-1}=s} (\delta_d\tens c^{-1}\la u)
 \delta_{ t , c^{-1}\ra(s\la u) }
\end{eqnarray*}
\begin{eqnarray*}
\big( c\ra(db\la u)  \big)^{-1} & = & \big( (cdb\ra u)(db\ra u)^{-1}
 \big)^{-1}  \\
& = &(db\ra u)(cdb\ra u)^{-1}\ =\ c^{-1}\ra (cdb\la u)\ =\
c^{-1}\ra (s\la u) \ .
\end{eqnarray*}
Finally to obtain the desired action we solve for $c$  which is fixed by the
delta function inside the summation. We have $t=c^{-1}\ra(s\la u)$, so
$c^{-1}=t\ra(s\la u)^{-1}=t'$. Also $d=c^{-1}sc=t's{t'}^{-1}$, so we have
\[  (t\tens\delta_v)\la(\delta_s\tens u)\ =\
\delta_{ v , (s\la u)^{-1}u } \big(\delta_{t's{t'}^{-1}}\tens t'\la u\big)\ .\]
\endproof

\proof :  The action of $H^*$ on $H$. First we calculate
$(\Delta\tens 1)\Delta(h)$ with $h=t\tens\delta_v\in H$:
\[
(\Delta\tens 1)\Delta(h)\ =\  \sum_{wzy=v} t\tens\delta_w\tens
t\ra w\tens\delta_z\tens  t\ra wz\tens\delta_y
\]
This gives summands $h_{(1)}=t\tens\delta_w$,
$h_{(2)}=t\ra w\tens\delta_z$,
and $h_{(3)}=t\ra wz\tens\delta_y$, which we substitute into the first
equation
of (3) to get
\begin{eqnarray*}
(t\tens\delta_v)\ra(\delta_s\tens u)& = & \sum_{wzy=v}(t\ra w\tens\delta_z)
\big<  S(t\tens\delta_w) (t\ra wz\tens\delta_y) , \delta_s\tens u \big>  \\
& = & \sum_{wzy=v}(t\ra w\tens\delta_z)
\big<  (t\ra w)^{-1}\tens\delta_{(t\la w)^{-1}}) (t\ra wz\tens\delta_y) ,
\delta_s\tens u \big>  \\
& = & \sum_{wzy=v}(t\ra w\tens\delta_z)
 \delta_{ (t\la w)^{-1}, (t\ra wz)\la y }\big<
(t\ra w)^{-1}(t\ra wz)\tens\delta_y  ,  \delta_s\tens u \big>  \\
& = & \sum_{wzy=v}(t\ra w\tens\delta_z)
 \delta_{ (t\la w)^{-1}, (t\ra wz)\la y }
\delta_{ (t\ra w)^{-1}(t\ra wz) , s }  \delta_{ y, u }  \\
& = & \sum_{wz=vu^{-1}}(t\ra w\tens\delta_z)
 \delta_{ (t\la w)^{-1}, (t\ra vu^{-1})\la u }
\delta_{ (t\ra w)^{-1}(t\ra vu^{-1}) , s }
 \end{eqnarray*}
Next we solve these equations for $w$ and $z$. We need the following
identities for double cross product groups:
\[
(t\la w)^{-1}\ = \ (t\ra w)\la w^{-1}  \quad,\quad
(t\ra w)^{-1}\ = \ t^{-1}\ra (t\la w)
 \]
Now if $(t\la w)^{-1}= (t\ra vu^{-1})\la u$ and
$s=(t\ra w)^{-1}(t\ra vu^{-1}) $, then $(t\ra w)\la w^{-1}= (t\ra vu^{-1})\la
u$, and so
\[
 w^{-1}\ =\  (t\ra w)^{-1}\la \big((t\ra vu^{-1})\la u\big)  \ =\ s\la u
 \]
Thus $w=(s\la u)^{-1}$, $z=(s\la u)vu^{-1}$, and $t\ra w=t\ra (s\la
u)^{-1}=t'$.
To simplify the second delta function in the summation we note that
\[   (t\ra w)^{-1}(t\ra vu^{-1})\ =\ t^{-1}(t\ra v)\ra u^{-1}\ .\]
This means that
\[
\delta_{ (t\ra w)^{-1}(t\ra vu^{-1}) , s }  \ = \
\delta_{s, t^{-1}(t\ra v)\ra u^{-1} }  \ =\
\delta_{s\ra u, t^{-1}(t\ra v) } \ =\
\delta_{t(s\ra u), t\ra v } \ ,
 \]
and this proves the formula for the action.
\endproof

\begin{propos} The quantum double $D(kM\cobicross k(G))$ is generated by
$kM\cobicross k(G)$ and $k(M)\bicross kG$ as sub-Hopf algebras with cross
relations defined by the product
\[ (1\tens t\tens\delta_v)(\delta_s\tens u\tens 1)=\delta_{t's(t\ra
vu^{-1})^{-1}}\tens (t\ra vu^{-1})\la u\tens t'\tens \delta_{(s\la u)vu^{-1}}\
,\]
where $t'=t\ra (s\la u)^{-1}$.
\end{propos}
\proof
Let $h=t\tens \delta_v\in H$, and $b=\delta_s\tens u\in H^*$. We calculate
$(1 \tens h)(b\tens 1)$ using (2). We have
\begin{eqnarray*}
\Delta(t\tens \delta_v)  & = &
\sum_{xy=v} t\tens \delta_x \tens t\ra x\tens \delta_y  \\
\Delta(\delta_s\tens u)  & = &
\sum_{wz=s} \delta_w \tens z\la u\tens \delta_z \tens u  \ ,
 \end{eqnarray*}
so that, by lemma 3.1,
\begin{eqnarray*}
h_{(1)}\la b_{(1)}  & = &
 (t\tens \delta_x) \la (  \delta_w \tens z\la u )  \\
 & = & \delta_{  z\la u  , (wz\la u)x   }.
 (\delta_{t'w{t'}^{-1}} \tens t'z\la u )  \\
 & = & \delta_{  z\la u  , (s\la u)x   }.
 (\delta_{t'w{t'}^{-1}} \tens t'z\la u )
 \end{eqnarray*}
where $t'=t\ra(s\la u)^{-1}$, since $wz=s$. Also
\[
h_{(2)}\ra b_{(2)}  \ = \
 \delta_{  t\ra xy  , (t\ra x) (z\ra u)   }.
 t'' \tens \delta_{(z\la u )yu^{-1} } \ =\
 \delta_{  t\ra v  , (t\ra x) (z\ra u)   }.
 t'' \tens \delta_{(z\la u )yu^{-1} }
 \]
where $t''=(t\ra x)\ra(z\la u)^{-1}$. Substituting into formula (2) gives
\[   (1\tens t\tens \delta_v) (\delta_s\tens u\tens 1) \ =\
\sum_{xy=v,\ wz=s}\delta_{z\la u,(s\la u)x}
\delta_{t\ra v,(t\ra x)(z\ra u)}
\big(\delta_{t'w{t'}^{-1}}\tens t'z\la u \tens t''
\tens \delta_{(z\la u)yu^{-1}} \big)
 \]
Next we solve
\[ z\la u\ =\ (s\la u)x\quad (a) \quad\text{and}\quad
 t\ra v\ = \ (t\ra x)(z\ra u)\quad (b)\ \]
 to find $t'z$, $t'w{t'}^{-1}$, $t''$ and $(z\la u)y$ in terms of $t$,
$v$, $s$ and $u$ alone. First we calculate
\begin{eqnarray*}
t'z\ra u& = & \big(t'\ra (z\la u)\big) \big(z\ra u\big)  \\
& = & \big((t\ra(s \la u)^{-1})\ra (z\la u)\big) \big(z\ra u\big) \\
& = & \big(t\ra(s \la u)^{-1}(z\la u)\big) \big(z\ra u\big) \\
& = & \big(t\ra x\big) \big(z\ra u\big)\ =\ t\ra v\ .
 \end{eqnarray*}
Thus $t'z=(t\ra v)\ra u^{-1}=t\ra vu^{-1}$. Also using (a) and $xy=v$ we have
\begin{eqnarray*}
 (z\la u)y& = &(s \la u)xy\ =\ (s \la u)v\ ,\ \text{and} \\
t''& = &(t\ra x)\ra (z\la u)^{-1} \ =\ t\ra x(z\la u)^{-1} \\
& = &t\ra (s\la u)^{-1}\ =\ t'\ .
 \end{eqnarray*}
Finally we have, using $wz=s$,
\[
t'w{t'}^{-1}\ =\ t'(wz)z^{-1}{t'}^{-1}\ =\ t's(t'z)^{-1}\ =\
t's(t\ra vu^{-1})^{-1}\ .
\]
 \endproof

We now give a relation between our order reversing group automorphisms on
$X$
 and the quantum double on $H=kM\cobicross k(G)$.
Remember that we had Hopf algebra isomorphisms
$\tilde\theta:H\to H^*$ and $\tilde\theta:H^*\to H$ given by
\[  \tilde\theta(s\tens\delta_u)\ =\ \delta_{\theta(s\la u)}\tens\theta(s\ra u)
\quad\text{and}\quad
\tilde\theta(\delta_s\tens u)\ =\ {\theta(s\la u)}\tens\delta_{\theta(s\ra u)}\
{}.
\]
In addition we have the order reversing map $\tau:H\tens H^*\to H^*\tens H$
 given by $\tau(h\tens b)=b\tens h$.

\begin{remark} We have already proved that $\tilde\theta:H\to H^*$
is a  Hopf algebra isomorphism. The given map $\tilde\theta:H^*\to H$ is the
inverse of $\widetilde{\theta^{-1}}:H\to H^*$, and so is also a
Hopf algebra isomorphism. We note that for $a\in H^*$ and $h\in H$,
\[
\big<\tilde\theta(a), \tilde\theta(h)\big>\ =\ \big<h, a\big>\ .
\]
\end{remark}
\proof
\begin{eqnarray*}
\big<\tilde\theta(\delta_{s}\tens u), \tilde\theta(t\tens \delta_v)\big>& = &
\big<\theta(s\la u)\tens \delta_{\theta(s\ra u)},
\delta_{\theta(t\la v)}\tens \theta(t\ra v)\big>  \\
& = & \delta_{\theta(s\la u),\theta(t\la v)}.
\delta_{\theta(s\ra u),\theta(t\ra v)} \ =\
\delta_{s\la u,t\la v}.\delta_{s\ra u,t\ra v} \\
& = & \delta_{s,t}.\delta_{u,v}\ =\
\big< t\tens \delta_v,\delta_{s}\tens u\big>\ .
\end{eqnarray*}
This equation has used the result that if $s\la u=t\la v$ and
$s\ra u=t\ra v$, then $su=(s\la u)(s\ra u)=(t\la v)(t\ra v)=tv$,
which means that $s=u$ and $t=v$.
\endproof

\begin{propos}
We have the following equations for $h\in H$ and $b\in H^*$:
\[
\tilde\theta(h\la b)\ =\ \tilde\theta b\ra \tilde\theta h \quad(a)
\quad\text{and}\quad
\tilde\theta(h\ra b)\ =\ \tilde\theta b\la \tilde\theta h \quad(b)\ .
\]
\end{propos}
\proof : (a) : We
 let $h=s\tens \delta_u$ and $b=\delta_t\tens v$, and find
\begin{eqnarray*}
\tilde\theta (\delta_t\tens v)\ra \tilde\theta (s\tens \delta_u)& = &
\delta_{v,(t\la v)u} . \big(\theta(t\la v)\ra \theta(s)^{-1} \big)
\tens \delta_{\theta\big(s(t\ra v)(s\ra u)^{-1}\big)}\ , \\
\tilde\theta\big( (s\tens \delta_u)\la  (\delta_t\tens v)\big)& = &
\delta_{v,(t\la v)u} . \theta(\bar s\la \bar u)
\tens \delta_{\theta(\bar s\ra \bar u)}\ ,
\end{eqnarray*}
where $\bar s=t't{t'}^{-1}$, $\bar u=t'\la v$, and $t'=s\ra(t\la v)^{-1}$. To
show that these are equal we calculate
\begin{eqnarray*}
\bar s\la \bar u & = & t't{t'}^{-1}\la(t'\la v)\ =\ t't\la v\ =\ t'\la(t\la v)
 \\ & = & \big(s\ra(t\la v)^{-1}\big)\la(t\la v)\ =\
\big(s\la(t\la v)^{-1}\big)^{-1}\ ,  \\
\theta(\bar s\la \bar u) & = & \theta\big(s\la(t\la v)^{-1}\big)^{-1}\ =\
\theta(t\la v)\ra\theta(s^{-1})\ .
\end{eqnarray*}
Also assuming that $v=(t\la v)u$, we have
\begin{eqnarray*}
\bar s\ra \bar u & = &
t't{t'}^{-1}\ra(t'\la v)\ =\  \big((t't{t'}^{-1})t'\ra v\big)(t'\ra
v)^{-1} \\
 & = & (t't\ra v)(t'\ra v)^{-1}\ =\
(t't\ra v)\big((s\ra(t\la v)^{-1})\ra v\big)^{-1} \\
 & = & (t't\ra v)\big((s\ra uv^{-1})\ra v\big)^{-1} \ =\
(t't\ra v)(s\ra u)^{-1} \\
 & = & (t'\ra (t\la v))(t\ra v)(s\ra u)^{-1}\ =\
\big((s\ra(t\la v)^{-1})\ra (t\la v)\big)(t\ra v)(s\ra u)^{-1}  \\
 & = & s(t\ra v)(s\ra u)^{-1}\ .
\end{eqnarray*}
\endproof
\proof : (b) :
We have
\[
\tilde\theta(\delta_t\tens v)\la \tilde\theta(s\tens \delta_u)\ =\
(\bar t\tens \delta_{\bar v})\la (\delta_{\bar s}\tens\bar u)\ =\
\delta_{\bar u,(\bar s\la \bar u)\bar v}.\delta_{{\bar t}'\bar
 s{\bar t}^{'-1}} \tens {\bar t}'\la \bar u
\ ,\]
where $\bar t=\theta(t\la v)$,
$\bar v=\theta(t\ra v)$,
$\bar s=\theta(s\la u)$,
$\bar u=\theta(s\ra u)$, and ${\bar t}'=\bar t\ra (\bar s\la \bar u)^{-1}$.
Then
\begin{eqnarray*}
\bar s\la \bar u& = & \theta(s\la u)\la \theta(s\ra u)\ =\ \theta(s) \ ,\\
{\bar t}'& = &\bar t\ra \theta(s)^{-1}  \ =\ \theta(t\la v)\ra \theta(s)^{-1}
\ =\ \theta\big(s\la(t\la v)^{-1}\big)^{-1}\ ,\\
\delta_{\bar u,(\bar s\la \bar u)\bar v}& = &
\delta_{\theta(s\ra u),\theta(s)\theta(t\ra v)}\ =\
\delta_{s\ra u,s(t\ra v)}\ .
\end{eqnarray*}Next we calculate
\begin{eqnarray*}
\tilde\theta\big((s\tens \delta_u)\ra(\delta_t\tens v)\big)& = &
\delta_{s\ra u,s(t\ra v)}
\tilde\theta\big( s\ra (t\la v)^{-1}\tens \delta_{(t\la v)uv^{-1}} \big) \\
& = & \delta_{s\ra u,s(t\ra v)}
\big( \delta_{\theta(m\la g)}\tens \theta(m\ra g) \big)\ ,
\end{eqnarray*}
where, if $s\ra u=s(t\ra v)$,
\begin{eqnarray*}
m\ra g & = & \big(s\ra (t\la v)^{-1}\big)\ra (t\la v) uv^{-1}\ =\
s\ra  uv^{-1}\ ,\\
m\la g & = & \big(s\ra (t\la v)^{-1}\big)\la (t\la v) uv^{-1}  \\
& = & \big( (s\ra (t\la v)^{-1}) \la (t\la v)  \big)
\big(\big(  (s\ra (t\la v)^{-1}) \ra (t\la v)  \big)\la uv^{-1}\big)   \\
& = & \big( s \la (t\la v)^{-1}  \big)^{-1}
\big(  (s\la uv^{-1}\big) \\
& = & \big( s \la (t\la v)^{-1}  \big)^{-1} (s\la u)
\big(  (s\ra u)\la v^{-1}\big) \\
& = & \big( s \la (t\la v)^{-1}  \big)^{-1} (s\la u)
\big( s\la\big( (t\ra v)\la v^{-1}\big)\big) \\
& = & \big( s \la (t\la v)^{-1}  \big)^{-1} (s\la u)
\big( s \la (t\la v)^{-1}  \big)\ .
\end{eqnarray*}
Now we have to show that ${\bar t}'\bar
 s{{\bar t}}^{'-1}=\theta(m\la g)$ (which is now automatic), and that
\begin{eqnarray*}
{\bar t}'\la \bar u& = &
\theta\big((s\ra u)\la v^{-1}\big)^{-1}\la\theta(s\ra u) \\
& = & \big(\theta(v)\ra \theta(s\ra u)^{-1}\big)\la\theta(s\ra u) \\
& = & \big(\theta(v)\la \theta(s\ra u)^{-1}\big)^{-1} \\
& = & \theta\big((s\ra u)\ra v^{-1}\big)^{-1}\ =\ \theta(s\ra uv^{-1})\ =\
\theta(m\ra g) \end{eqnarray*}
\endproof

\begin{remark}
For $h\in H$ and $b\in H^*$, we have $\big<S(h),S(a)\big>=\big<h,a\big>$.
This also
implies that $S^2$ is the identity on $H$ and $H^*$.
\end{remark}

\begin{propos} The map $\psi=\tau(\tilde\theta \tens \tilde\theta):
D(H)\to D(H)$ is an anti-algebra isomorphism, a
 coalgebra isomorphism, and preserves the unit, the counit, and the anitpode.
\end{propos}
\proof
First we show that $\psi$ is an anti-algebra map. On the one hand we have,
using
proposition 3.4,
\begin{eqnarray*}
\psi\big((a\tens h)(b\tens g)\big) & = & \psi \Big(
\sum (h_{(1)} \la b_{(1)})a\tens (h_{(2)} \ra b_{(2)})g \Big)  \\
& = & \sum
\tilde\theta(h_{(2)} \ra b_{(2)})\tilde\theta(g)
\tens
\tilde\theta(h_{(1)} \la b_{(1)})\tilde\theta(a)
  \\
& = & \sum \big(\tilde\theta g\tens 1\big)  \big(
\tilde\theta(h_{(2)} \ra b_{(2)})
\tens
\tilde\theta(h_{(1)} \la b_{(1)})
\big)   \big(1\tens\tilde\theta a \big)  \\
& = & \sum \big(\tilde\theta g\tens 1\big)  \big(
(\tilde\theta b)_{(2)} \la (\tilde\theta h)_{(2)}
\tens
(\tilde\theta b)_{(1)} \ra (\tilde\theta h)_{(1)}
\big)   \big(1\tens\tilde\theta a \big)
\end{eqnarray*}
On the other hand,
\begin{eqnarray*}
\psi(b\tens g)\psi(a\tens h) & = & (\tilde\theta g\tens \tilde\theta b)
 (\tilde\theta h\tens \tilde\theta a)  \\
& = & \sum
\big((\tilde\theta b)_{(1)} \la (\tilde\theta h)_{(1)}\big)
\tilde\theta g\tens
\big((\tilde\theta b)_{(2)} \ra (\tilde\theta h)_{(2)}\big)
\tilde\theta a \\
& = & \sum \big(\tilde\theta g\tens 1\big)  \big(
(\tilde\theta b)_{(1)} \la (\tilde\theta h)_{(1)} \big)
\tens
\big((\tilde\theta b)_{(2)} \ra (\tilde\theta h)_{(2)}
\big)   \big(1\tens\tilde\theta a \big)\ .
\end{eqnarray*}
To complete the calculation we note the identity [5, Chapter 7]
\[
\sum h_{(1)}\la a_{(1)} \tens h_{(2)}\ra a_{(2)} \ =\
\sum h_{(2)}\la a_{(2)} \tens h_{(1)}\ra a_{(1)}\ .
\]

We now check the condition for $\psi$ to be a coalgebra map, i.e.\
$(\psi\tens\psi)\Delta=\Delta \psi$.
\begin{eqnarray*}
\Delta \psi(a\tens h) & = & \Delta (\tilde\theta h\tens \tilde\theta a) \\
& = &\sum (\tilde\theta h)_{(1)} \tens
(\tilde\theta a)_{(1)} \tens
(\tilde\theta h)_{(2)} \tens
(\tilde\theta a)_{(2)}  \\
& = &\sum \tilde\theta (h_{(1)}) \tens
\tilde\theta (a_{(1)}) \tens
\tilde\theta (h_{(2)}) \tens
\tilde\theta (a_{(2)})  \\
& = &\sum \psi\big( a_{(1)} \tens h_{(1)} \big) \tens
\psi\big(  a_{(2)} \tens h_{(2)} \big) \\
& = & (\psi\tens \psi)\sum  a_{(1)} \tens h_{(1)}  \tens
  a_{(2)} \tens h_{(2)}  \\
& = & (\psi\tens \psi)\Delta  (a \tens h) \ .
\end{eqnarray*}
To check that the antipode is preserved,
\begin{eqnarray*}
S \psi(a\tens h) & = & S(\tilde\theta h\tens \tilde\theta a) \\
& = & (1\tens S\tilde\theta a)( S\tilde\theta h \tens 1) \\
& = & (1\tens \tilde\theta S a)( \tilde\theta S h \tens 1) \\
& = & \psi( \tilde\theta S a\tens 1)\psi(1\tens \tilde\theta S h) \\
& = & \psi\big( (1\tens \tilde\theta S h)( \tilde\theta S a\tens 1)\big) \\
& = & \psi S(a\tens h)\ .
\end{eqnarray*}

Finally we check the effect of $\psi$ on the co-unit and unit of $D(H)$, using
the fact that $\tilde\theta$ is a Hopf algebra isomorphism:
\begin{eqnarray*}
\eps_{D(H)}\psi(a\tens h) & = & \eps_{D(H)}\big(\tilde\theta h \tens
\tilde\theta a \big) \\
& = &(\eps \tilde\theta h )(\eps \tilde\theta a ) \ =\
(\eps  h )(\eps  a )\\
& = &
\eps_{D(H)}\big(a \tens h \big)\ , \\
\psi(1\tens 1) & = & \tilde\theta(1)\tens \tilde\theta(1)\ =\ 1\tens 1\ .
\end{eqnarray*}
\endproof

\begin{example} In the case of our $\Z_{2n}\Z_{2n}$ example,
we can write out the
actions as (we use additive operations in $Z_{2n}$)
\begin{eqnarray*}
(t\tens\delta_v)\la (\delta_s\tens u)& =&
\left\{ \begin{array}{lll}
\delta_{0,v}.\delta_s\tens u  & t\ \text{ even}  & s\ \text{ even} \\
\delta_{0,v}.\delta_s\tens (-u)  & t\ \text{ odd}  & s\ \text{ even} \\
\delta_{2u,v}.\delta_s\tens u  & t\ \text{ even}  & s\ \text{ odd} \\
\delta_{2u,v}.\delta_s\tens (-u)  & t\ \text{ odd}  & s\ \text{ odd}
\end{array} \right.
\\
(t\tens\delta_v)\ra (\delta_s\tens u)& =&
\left\{ \begin{array}{lll}
\delta_{0,s}.t\tens \delta_{v}  & u\ \text{ even}  & v\ \text{ even} \\
\delta_{0,s}.(-t)\tens \delta_{v}  & u\ \text{ odd}  & v\ \text{ even} \\
\delta_{-2t,s}.t\tens \delta_{v}  & u\ \text{ even}  & v\ \text{ odd} \\
\delta_{2t,s}.(-t)\tens \delta_{v}  & u\ \text{ odd}  & v\ \text{ odd}
\end{array} \right.
\end{eqnarray*}
\end{example}

\section{Representations and Braiding}

The representations of the quantum double of any Hopf algebra form a braided
tensor category, and
for this reason the category is especially interesting. By the
double cross product
construction we see that we have such a braided category for every group
factorisation. Many results and applications flow from this. We shall mention
only one or two of them.

The respresentations of $D(H)$ are evident from its description as a double
cross product. namely, they are left $H$ and right $H^*$-modules $W$ which
are
compatible in such a way that
\eqn{bimod}{ h\la (w\ra a)=\sum( (h\o\ra a\o)\la w) \ra(h\t\la a\t)}
which can be further computed in terms of the mutual coadjoint actions. We
freely identify a left $H^{*\rm op}$ module as a right $H^*$-module. The
braiding among such modules $V,W$ is
\eqn{braiding}{\Psi_{V,W}(v\tens w)=\sum _a e_a\la w\tens v\ra f^a}
where $\{e_a\}$ is a basis of $H$ and $\{f^a\}$ is a dual basis. The
corresponding action of $D(H)$ is $(a\tens h)\la w=(h\la w)\ra a$.

Finally, a canonical representation of $D(H)$ (aside from the obvious
left-regular representation on itself) is on $H$. This is motivated by
thinking
of $D(H)$ as a semidirect product of the type arising in quantum mechanics
(in
some cases it is) and can be called the `Schroedinger representation' of the
quantum double. Explicitly, it is
\eqn{canrep}{ h\la g=\sum h\o g Sh\t,\quad h\ra a=\sum \<a, h\o\>h\t,\qquad
\forall h,g\in H,\ a\in H^*}
which are respectively the quantum adjoint action of $H$ on itself and the
right coregular action of $H^*$ on $H$ induced by the left regular coaction.
The braiding in this case reduces to
\eqn{canbraid}{ \Psi_{H,H}(h\tens g)=\sum h\o g Sh\t\tens h\th.}
This is a canonical braiding (i.e. a solution of the celebrated quantum
Yang-Baxter equations) associated to any Hopf algebra. For this last
observation the Hopf algebra need not be finite-dimensional.

\begin{propos} The representations of the double $D(kM\cobicross k(G))$
are in
one-one correspondence with vector spaces $W$ which are

(i) $G$-graded left $M$-modules such that $|t\la w|=t\la |w|$ for all $t\in M$,
where $|\ |$ denotes the $G$-degree of a homogeneous element $w\in W$.

(ii)  $M$-graded right $G$-modules such that $\<w\ra u\>=\<w\>\ra u$ for all
$u\in G$, where $\<\ \>$ denotes the $M$-degree of a homogeneous element
$w\in
W$.

(iii) mutually `cross modules' according to
\[ \<t\la w\>=t\<w\>(t\ra |w|)^{-1},\quad |w\ra u|=(\<w\>\la u)^{-1}|w|u\]
on homogeneous elements.

(iv) $G-M$-`bimodules' according to
\[ (t\ra (\<w\>\la u))\la(w\ra u)=(t\la w)\ra ((t\ra|w|)\la u)\]

The corresponding action of the double and the induced braiding are
\[ (t\tens \delta_v)\la w=t\la w\delta_{v,|w|},\quad w\ra (\delta_s\tens
u)=\delta_{s,\<w\>} w\ra u.\]
\[ \Psi_{V,W}(v\tens w)=\<v\>\la w\tens v\ra |w|.\]
\end{propos}
\proof (if direction)
 Suppose that $W$ is a vector space satisfying conditions {\it (i)} to
{\it (iv)}. First we show that {\it (i)} ensures that $W$ is a left
$H$-module. Since $W$ is $G$-graded, we have $W=\oplus_{g\in G} W_g$,
where
$|w|=g$ if and only if $w\in W_g$. Also $W$ is a left $k(G)$-module by
$\delta_g\la w=\delta_{g,|w|}w$ for homogenous $w$ (i.e.\ $w\in W_{|w|}$).
Additionally from {\it (i)}, $W$ is a left $M$-module, and hence a left
$kM$-module. Now define an action of $H=kM\cobicross k(G)$ on $W$ by
\[
(t\tens \delta_v)\la w\ =\ t\la(\delta_v\la w)\ =\ \delta_{v,|w|}t\la w
\ ,\quad \text{where }w\in W_{|w|}\quad (a).
\]
To see that this is a left action, calculate
\begin{eqnarray*}
\big((s\tens \delta_u)(t\tens \delta_v)\big)\la w & = & \delta_{u,t\la v}
(st\tens \delta_v)\la w   \\
& = & \delta_{u,t\la v} \delta_{v,|w|} st\la w \ ,\\
(s\tens \delta_u)\la \big((t\tens \delta_v)\la w\big) & = &
\delta_{v,|w|}(s\tens \delta_u)\la (t\la w)  \\
& = &  \delta_{v,|w|} \delta_{u,|t\la w|} st\la w \ .
\end{eqnarray*}
These are equal since $|t\la w|=t\la |w|$. Also
\[
1_H\la w\ =\ \Big(\sum_u e\tens\delta_u\Big)\la u\ = \ \sum_u
 \delta_{u,|w|} e\la w\ =\  \delta_{u,|w|} w\ =\ w\ .
\]
Next we show that {\it (ii)} ensures that $W$ is a right $H^*$-module.
Since $W$ is $M$-graded, $W=\oplus_{s\in M} W_s$, where
$<w>=s$ if and only if $w\in W_s$.
 Also $W$ is a right $k(M)$-module by
$w\ra \delta_s=\delta_{s,<w>}w$ for homogenous $w$ (i.e.\ $w\in W_{<w>}$).
By assumption $W$ is also a right $kG$-module. Now define a right action of
$H^*=k(M)\bicross kG$ on $W$ by \[
w\ra (\delta_s\tens u)\ =\ (w\ra \delta_s)\ra u\ =\ \delta_{s,<w>} w\ra u
\ ,\quad \text{where }\ w\in W_{<w>}\quad (b).
\]
To see that this is a right action, calculate
\begin{eqnarray*}
w\ra \big((\delta_s\tens u)(\delta_t\tens v)\big) & = &
\delta_{s\ra u,t} w\ra (\delta_s\tens uv)  \\
& = &\delta_{s\ra u,t}  \delta_{s,<w>} w\ra uv  \ ,\\
\big(w\ra (\delta_s\tens u)\big)\ra (\delta_t\tens v)& = &
 \delta_{s,<w>} (w\ra u)\ra  (\delta_t\tens v) \\
& = &  \delta_{s,<w>} \delta_{t,<w\ra u>} w\ra uv\ .
\end{eqnarray*}
These are equal since $<w\ra u>=<w>\ra u$.

The action of $D(H)$ on $W$ is defined by $(a\tens h)\la w=(h\la w)\ra a$.
Now we must check that this does define a left action of $D(H)$. If
$(1\tens h)(a\tens 1)=a'\tens h'$ in $D(H)$, then we show that
$(1\tens h)\la\big((a\tens 1)\la w\big)=(a'\tens h')\la w$, i.e.\ that
$h\la(w\ra a)=(h'\la w)\ra a'$. Putting $h=t\tens \delta_v$
and $a=\delta_s\tens u$, assuming that $w$ is homogenous,
 and using formulae (a) and (b), we get the equivalent
 conditions
\[
 \delta_{s,<w>} \delta_{v,|w\ra u|} t\la(w\ra u)   \ =\
  \delta_{v',|w|} \delta_{s',<t'\la w>} (t'\la w) \ra u'
\quad (c)\ ,
\]
where $(1\tens t\tens \delta_v)(\delta_s\tens u\tens 1)=\delta_{s'}\tens u'
 \tens t'\tens \delta_{v'}$, and using Proposition 3.2,
\[
s'=\big(t\ra (s\la u)^{-1}\big)s(t\ra vu^{-1})^{-1},\quad u'=
(t\ra vu^{-1})\la u,\quad t'=t\ra (s\la u)^{-1},\quad v'=(s\la u)vu^{-1}\ .
\]
It will be useful to calculate
\begin{eqnarray*}
t'\ra v' & = & t\ra (s\la u)^{-1}\ra (s\la u)vu^{-1}\ =\ t\ra vu^{-1}\ ,\\
s'\la u' & = & \big(t\ra (s\la u)^{-1}\big)s(t\ra vu^{-1})^{-1} \la
\big((t\ra vu^{-1})\la u\big)  \\
& = & \big(t\ra (s\la u)^{-1}\big)s \la u \ =\
 \big(t\ra (s\la u)^{-1}\big)\la(s \la u)  \\
& = & \big(t\la (s\la u)^{-1}\big)^{-1} \ ,\\
t'\la v' & = &  \big(t\ra (s\la u)^{-1}\big)\la\big((s\la u)vu^{-1}\big)
\ =\  \big(\big(t\ra (s\la u)^{-1}\big)\la (s\la u)\big)
\big(\big(t\ra (s\la u)^{-1}\big)\ra (s\la u)\big)
\la vu^{-1}\big) \\
 & = & \big(t\ra (s\la u)^{-1}\big)^{-1}t\la  vu^{-1}\ =\ (s'\la u')
t\la  vu^{-1}\ , \\
s'\ra u'& = &\big(t\ra (s\la u)^{-1}\big)s(t\ra vu^{-1})^{-1}    \ra
\big((t\ra vu^{-1})\la u\big)  \\
& = & \big(  \big(t\ra (s\la u)^{-1}\big)s \ra u \big)
 \big( (t\ra vu^{-1})^{-1} \ra  \big( (t\ra vu^{-1}) \la u \big) \big) \\ & = &
 \big(t\ra (s\la u)^{-1}\ra (s\la u) \big) (s\ra u)(t\ra vu^{-1}\ra u)^{-1} \\
& = & t(s\ra u)(t\ra v)^{-1} \ .
\end{eqnarray*}
Now we use {\it (iii)} and {\it (iv)} to show that the condition (c) holds.
First we prove that
\[
s=<w> \text{ and } v=|w\ra u| \quad \text{ if\ and\ only\ if } \quad
v'=|w| \text{ and } s'=<t'\la w>\ .
\]
If $s=<w>$ and $v=|w\ra u|$, then
\begin{eqnarray*}
v' & = & (s\la u)vu^{-1} \ =\ (<w>\la u)vu^{-1} \\
& = & |w|u|w\ra u|^{-1}vu^{-1} \quad \text{by {\it (iii)}}  \\
& = & |w|uv^{-1}vu^{-1}\ =\ |w|\ ,\\
s' & = & t's(t\ra vu^{-1})^{-1}\ =\ t's(t'\ra v')^{-1}
\ =\ t's(t'\ra |w|)^{-1}\ =\ t'<w>(t'\ra |w|)^{-1} \\
& = & <t'\la w> \quad \text{by {\it (iii)}} \ .
\end{eqnarray*}
Conversely, if $v'=|w|$ and $s'=<t'\la w>$, then
\begin{eqnarray*}
s & = & t^{'-1}s'(t\ra vu^{-1})\ =\  t^{'-1}<t'\la w>(t\ra vu^{-1})  \\
& = &  t^{'-1}t'< w>(t'\ra|w|)^{-1}(t\ra vu^{-1})
 \quad \text{by {\it (iii)}}  \\
& = & < w>(t'\ra v')^{-1}(t'\ra v') \ =\ <w>\ ,\\
v & = & (s\la u)^{-1}v'u \ =\ (<w>\la u)^{-1}|w|u  \\
& = & |w\ra u| \quad \text{by {\it (iii)}} \ .
\end{eqnarray*}
Now assuming that $s=<w>$ and $v=|w\ra u|$, we show that
$t\la (w\ra u)=(t'\la w)\ra u'$. If we define
 $\bar t=t\ra(<w>\la u)^{-1}$, then
\begin{eqnarray*}
t\la (w\ra u)& = & \big(\bar t \ra  (<w>\la u)\big) \la (w\ra u)  \\
& = & (\bar t\la w)\ra \big((\bar t\ra |w|)\la u\big)
\quad \text{by {\it (iv)}}  \\
& = & (\bar t\la w)\ra \big(((t\ra(<w>\la u)^{-1})\ra |w|)\la u\big)  \ =\
(\bar t\la w)\ra \big((t\ra(<w>\la u)^{-1}|w|)\la u\big)  \\
& = & (\bar t\la w)\ra \big((t\ra |w\ra u|u^{-1})\la u\big)
\quad \text{by {\it (iii)}}  \\
& = & (\bar t\la w)\ra \big((t\ra vu^{-1})\la u\big)\ =\
(\bar t\la w)\ra u'  \\
& = & ((t\ra(<w>\la u)^{-1})\la w)\ra u'\ =\
((t\ra |w\ra u|u^{-1}|w|^{-1})\la w)\ra u'  \quad \text{by {\it (iii)}}  \\
& = & ((t\ra vu^{-1}v^{'-1})\la w)\ra u' \ =\
((t\ra vu^{-1}uv^{-1}(s\la u)^{-1})\la w)\ra u'  \\
& = &  ((t\ra (s\la u)^{-1})\la w)\ra u' \ =\ (t'\la w)\ra u'\quad (d)\ .
\end{eqnarray*}
\endproof
\proof (only if direction):  Suppose now that we have a representation of
$D(H)$ on a vector space $W$. Then $W$ is a left $H$-module. But as an
algebra
$H=kM\cobicross k(G)$, the semidirect product of $kM$ and $k(G)$. This
means
that $W$ is a left $M$-covariant left $k(G)$-module. But a  $k(G)$-module
is precisely a $G$-graded vector space, where the action of $\delta_g$ is
$\delta_g\la w=\delta_{g,|w|}w$, for $w$ homogenous. This and
$M$-covariance
gives {\it (i)}. The action of $H$ on $W$ is $(t\tens\delta_v)\la w =
\delta_{v,|w|}t\la w$, where $|.|$ is the $G$-grading on $W$.

Similarly $W$ is a right $H^*$-module, and $H^*=k(M)\bicross kG$. By the
above
reasoning, $W$ is a right $G$-covariant right $k(M)$-module. This is
condition
 {\it (ii)}. The action of $H^*$ on $W$ is $w\ra (\delta_s\tens u)=
\delta_{s,<w>}w\ra u$, where $<.>$ is the $M$-grading on $W$.

Because we have a representation of $D(H)$ on $W$, and have determined
 the corresponding actions of $H$ and $H^*$, we know that for
homogenous $w$,
\[
 \delta_{s,<w>} \delta_{v,|w\ra u|} t\la(w\ra u)   \ =\
  \delta_{v',|w|} \delta_{s',<t'\la w>} (t'\la w) \ra u'
\quad (c)\ .
\]
We show that this implies  {\it (iii)} and  {\it (iv)}. Now assume that
$s=<w>$ and $v=|w\ra u|$. Then
\begin{eqnarray*}
|w|& = & v'\ =\ (s\la u)vu^{-1}\   \\
& = & (<w>\la  u)vu^{-1}\ =\  (<w>\la  u)|w\ra u|u^{-1}\ , \\
\text{so that}\quad |w\ra u| & = & (<w>\la  u)^{-1}|w|u\ ,\\
<w>& = & s\ =\ t^{'-1}s'(t\ra vu^{-1})\ =\  t^{'-1}<t'\la w>(t\ra vu^{-1})  \\
& = &  t^{'-1}<t'\la w>(t'\ra v')\ =\  t^{'-1}<t'\la w>(t'\ra |w|)\ , \\
\text{so that}\quad <t'\la w> & = & t'<w>(t'\ra |w|)^{-1}\ ,
\end{eqnarray*}
which gives  {\it (iii)}.

To prove  {\it (iv)}, note that by following the calculation (d) backwards
 to the second step (which is legitimate, since we now know that
 {\it (iii)} holds) gives
\[
(t'\la w)\ra u'\ =\ (\bar t\la w)\ra \big((\bar t\ra |w|)\la u\big)
\]
where $\bar t=t\ra(<w>\la u)^{-1}$. Then
\[
(t'\la w)\ra u'\ =\ t\la(w\ra u)\ =\ \big(\bar t\ra(<w>\la u)\big)\la (w\ra u)
\ ,\]
so that we have the following equation, which in turn gives {\it (iv)}:
\[
(\bar t\la w)\ra \big((\bar t\ra |w|)\la u\big)\ =\
\big(\bar t\ra(<w>\la u)\big)\la (w\ra u)\ .
\]

Let $V$ and $W$ be representations of $D(H)$. Then the induced braiding is
given by (5), i.e.
\begin{eqnarray*}
\Psi_{V,W}(v\tens w)& = & \sum(s\tens\delta_u)\la w
\tens v\ra (\delta_s\tens u)\ =\
\sum(\delta_{u,|w|} s\la w)
\tens (\delta_{s,<v>} v\ra  u) \\
& = & \<v\>\la w\tens v\ra |w|\ .
\end{eqnarray*}
 \endproof

Note that conditions (iii) are each a generalisation of Whitehead's notion
of crossed modules\cite{Whi:com},  but now the `crossed $G$-module' and the
`crossed $M$-module' are coupled unless one of the factors is trivial. For
example, if $M$ is trivial then we have exactly a right crossed $G$-module.
The
connection between quantum doubles of a single group and crossed modules
was
introduced in \cite{Ma:dou}, and the last proposition extends this point of
view to the quantum double of a
bicrossproduct associated to a matched pair of
groups. In view of this, we call modules obeying the conditions in the
proposition {\em bicrossed bimodules}.

\begin{propos} The Schroedinger representation of the quantum double on
$W=kM\cobicross k(G)$ and the induced braiding are
\[ (s\tens \delta_u)\la (t\tens\delta_v)=st(s\ra u)^{-1}\tens\delta_{(s\ra
u)\la v}\delta_{uv,t\la v},\quad  (t\tens\delta_v)\ra(\delta_s\tens
u)=\delta_{s,t} t\ra u\tens \delta_{u^{-1}v}.\]
\[ \Psi(s\tens\delta_u\tens t\tens\delta_v)=st {s'}^{-1}\tens\delta_{s'\la
v}\tens s'\tens \delta_{v(t\la v)^{-1}u};\quad s'=s\ra(t\la v)v^{-1}.\]
\end{propos}

\begin{corol} The gradings and $M-G$ actions for the Schroedinger
representation $W=kM\bicross k(G)$ according to Proposition~4.1 are
\[ |t\tens\delta_v|=(t\la v)v^{-1},\quad \<t\tens\delta_v\>=t\]
\[ s\la(t\tens\delta_v)=st{s'}^{-1}
\tens \delta_{s'\la v},\quad (t\tens\delta_v)\ra
u=t\ra u\tens \delta_{u^{-1}v};\quad s'=s\ra(t\la v)v^{-1}.\]
\end{corol}
\proof of proposition and corollary: First we use formula (6) to calculate the
action of $H$ on itself, taking $h=s\tens\delta_u$, $g=t\tens\delta_v$.
Since $\Delta(h)=\sum_{xy=u} s\tens\delta_x\tens s\ra x\tens\delta_y$,
\begin{eqnarray*}
(s\tens\delta_u)\la(t\tens\delta_v)& = & \sum_{xy=u}
(s\tens\delta_x)(t\tens\delta_v)S(s\ra x\tens \delta_y) \\
& = & \sum_{xy=u} (s\tens\delta_x)(t\tens\delta_v) \big(
(s\ra xy)^{-1}\tens\delta_{((s\ra x)\la y)^{-1}}
\big)  \\
& = & \sum_{xy=u} \delta_{x,t\la v} (st\tens\delta_v)   \big(
(s\ra u)^{-1}\tens\delta_{(s\la u)^{-1}(s\la x)}
\big)  \\
& = &  (st\tens\delta_v)   \big(
(s\ra u)^{-1}\tens\delta_{(s\la u)^{-1}(st\la v)}
\big)  \\
& = &  \delta_{v,(s\ra u)^{-1}\la (s\la u)^{-1}(st\la v) }   \big(
st(s\ra u)^{-1}  \tens  \delta_{(s\la u)^{-1}(st\la v)}
\big)  \ ,
\end{eqnarray*}
where we have used
\[
(s\ra x)\la y\ =\ (s\ra x)\la x^{-1}u\ =\ (s\la x)^{-1}(s\la u)\ .
\]
To simplify further, we calculate
\begin{eqnarray*}
(s\ra u)^{-1}\la (s\la u)^{-1}(st\la v) & = &
\big((s\ra u)^{-1}\la (s\la u)^{-1}\big)
\big(\big((s\ra u)^{-1}\ra (s\la u)^{-1}\big)\la (st\la v)\big) \\
& = & u^{-1}\big(s^{-1}\la (st\la v)\big)\ =\ u^{-1}(t\la v)\ ,
\end{eqnarray*}
so that we may assume that $t\la v=uv$, and then
\[
 (s\la u)^{-1}(st\la v) \ = \  (s\la u)^{-1}(s\la uv) \ =\ (s\la u)^{-1}
(s\la u)((s\ra u)\la v)\ =\ (s\ra u)\la v\ .
\]
Hence we obtain
\[
(s\tens\delta_u)\la(t\tens\delta_v)\ =\
 \delta_{uv,t\la v }   \big(
st(s\ra u)^{-1}  \tens  \delta_{(s\ra u)\la v}
\big)\ .
\]
Next we calculate the right coregular action of $H^*$ on $H$, using formula
(6) with $h=t\tens\delta_v$ and $a=\delta_s\tens u$,
\begin{eqnarray*}
(t\tens\delta_v)\ra (\delta_s\tens u)& =& \sum_{xy=v}
 \big<\delta_s\tens u,t\tens\delta_x\big>(t\ra x)\tens\delta_y  \\
& =& \sum_{xy=v}
 \delta_{s,t}\delta_{u,x}(t\ra x)\tens\delta_y\ =\  \delta_{s,t}
(t\ra u)\tens\delta_{u^{-1}v}\ .
\end{eqnarray*}

We calculate the induced braiding using the formula
\[ \Psi_{H,H}(v\tens w)\ =\ <v>\la w \tens v\ra |w|\ ,\]
from proposition 4.1. We could also use (7).
 But first we must determine the $M$ and $G$ gradings,
 and the actions of $M$ and $G$ on $H$. From the first part of the proof,
and proposition 4.1,
\begin{eqnarray*}
(s\tens \delta_u)\la(t\tens \delta_v) & = & \delta_{uv,t\la v}st(s\ra u)^{-1}
\tens \delta_{(s\ra u)\la v}  \ ,\\
(s\tens \delta_u)\la(t\tens \delta_v) & = & \delta_{u,|t\tens \delta_v|}
s\la (t\tens \delta_v)  \ ,
\end{eqnarray*}
so that
\begin{eqnarray*}
|t\tens \delta_v|& = & (t\la v)v^{-1}\ , \\
s\la (t\tens \delta_v)& = & st(s\ra u)^{-1}
\tens \delta_{(s\ra u)\la v} \\
& = & sts^{'-1}\tens \delta_{s'\la v}\ ,
\end{eqnarray*}
where $s'=s\ra u=s\ra(t\la v)v^{-1}$ (in view of the delta function,
we may take $u=(t\la v)v^{-1}$). Similarly we showed
\begin{eqnarray*}
(t\tens \delta_v)\ra (\delta_s\tens u) & = & \delta_{s,t}t\ra u
\tens \delta_{u^{-1} v}  \ ,\\
(t\tens \delta_v)\ra (\delta_s\tens u) & = & \delta_{s,<t\tens \delta_v>}
(t\tens \delta_v)\ra u  \ ,
\end{eqnarray*}
so that
\begin{eqnarray*}
<t\tens \delta_v>& = & t\ , \\
(t\tens \delta_v)\ra u& = & t\ra u
\tens \delta_{u^{-1} v}\ .
\end{eqnarray*}
Then the braiding is given by
\begin{eqnarray*}
\Psi_{H,H}(s\tens \delta_u\tens t\tens \delta_v)& = &
<s\tens \delta_u>\la (t\tens \delta_v)
\tens (s\tens \delta_u)\ra |t\tens \delta_v|  \\
& = &
s\la (t\tens \delta_v)
\tens (s\tens \delta_u)\ra (t\la v)v^{-1}  \\
& = &
sts^{'-1}\tens \delta_{s'\la v}
\tens s\ra(t\la v) v^{-1}\tens  \delta_{v(t\la v)^{-1}u}  \\
& = &
sts^{'-1}\tens \delta_{s'\la v}
\tens s'\tens  \delta_{v(t\la v)^{-1}u}
\end{eqnarray*}
 \endproof

\begin{example} For our $\Z_{2n}\Z_{2n}$ example the brading matrix
$\Psi$ is
\[
\Psi\big(s\tens\delta_u\tens t\tens\delta_v\big)\  = \
\left\{ \begin{array}{lll}
t\tens\delta_{v}\tens s\tens\delta_{u}
& t\ \text{ even}  & s\ \text{ even} \\
t\tens\delta_{v}\tens s\tens\delta_{u+2v}
& t\ \text{ odd}  & s\ \text{ even} \\
t\tens\delta_{-v}\tens s\tens\delta_{u}
& t\ \text{ even}  & s\ \text{ odd} \\
t\tens\delta_{-v}\tens s\tens\delta_{u+2v}
& t\ \text{ odd}  & s\ \text{ odd}
\end{array} \right.
\]
In the case $\Z_{6}\Z_{6}$ the minimal polynomial of $\psi$ is
\[ -1 - \lambda^2 - \lambda^4 + \lambda^8 + \lambda^{10} + \lambda^{12}\ .\]
\end{example}
 \proof
Remember that both actions of $\Z_{2n}$ on $\Z_{2n}$ are given by
an even element having trivial action, and an odd element reversing sign.
Then the gradings are given by the formulae $|t\tens \delta_v|=0$ if $t$
is even, and $|t\tens \delta_v|=-2v$ if $t$
is odd. Also $<t\tens \delta_v>=t$. This gives the formula for $\Psi$.

The effect of the
linear transformation is to permute the elements of the bases. If we have
 a vector space of dimension $m$ and a linear transformation which
cyclically permutes the basis vectors, the characteristic polynomial is
$\lambda^m-1$. Since all the roots are distinct, this is also the minimal
polynomial. If we have a vector space which is a direct sum of subspaces
where a linear transformation acts on each subspace, the minimal polynomial
is the least common multiple of the minimal polynomials on each
subspace.

To get the minimal polynomial of $\Psi$, we decompose the permutation of the
basis vectors into disjoint cycles. A cycle of length $m$ contributes
$\lambda^m-1$, and we take the least common multiple of the $\lambda^m-1$
polynomials over the cycles. The order of $\Psi$ is the
least common multiple of the lengths of the cycles.
Note that $\Psi$ preserves $s$ and $t$ in
$s\tens\delta_u\tens
t\tens\delta_v$,  and we can consider the different cases
even and odd for $s$ and $t$
separately.

The calculation of the sizes of the cycles is in general a non-trivial
problem in number theory. We shall only do the calculation for the case
$n=3$.

If both $s$ and $t$ are even, these subspaces are 1 dimensional (ie. no
effect). The polynomial is $\lambda-1$.

If $s$ is odd and $t$ even, we get the shift $(u,v) \mapsto
(-v,u)$, which has order 4, and subspaces of order 1, 2 and 4.
 The polynomial is $\lambda^4-1$.

If $s$ is even and $t$ odd, we get the shift $(u,v) \mapsto
(v,u+2v)$, which has order 8, and subspaces of order 1, 2 and 8.
 The polynomial is $\lambda^8-1$.

If both $s$  and $t$ are odd, we get the shift $(u,v) \mapsto
(-v,u+2v)$, which has order 6, and subspaces of order 1, 2, 3 and 6.
 The polynomial is $\lambda^6-1$.

 Taking the least common multiple of the $\lambda^m-1$ factors, we see that
$\Psi$ has minimal polynomial $(\lambda^8-1)(\lambda^6-1)/(\lambda^2-1)=
-1 - \lambda^2 - \lambda^4 + \lambda^8 + \lambda^{10} + \lambda^{12}$.
 \endproof

\section{$D(H)$ as a twisting of $D(X)$.}

We recall that $X$ is the double cross product group $GM$. This group has a
quantum double  $D(X)=k(X)\lcross kX$ which is an ordinary crossproduct,
It has operations
\[ (\delta_x\tens y)(\delta_a\tens b)=\delta_{y^{-1}xy,a}(\delta_x\tens yb),
\quad
\Delta(\delta_x\tens y)=\sum_{ab=x}\delta_a\tens y\tens \delta_b\tens y\]
\[ 1=\sum_x\delta_x\tens e,\quad \eps(\delta_x\tens y)=\delta_{x,e},\quad
S(\delta_x\tens y)=\delta_{y^{-1}x^{-1}y}\tens y^{-1},\quad
(\delta_x\tens y)^*=\delta_{y^{-1}x y}\tens y^{-1}
,\]
\[
R\ =\ \sum_{y,z} \delta_y\tens e \tens \delta_z\tens y\ .
\]

The representations of $D(X)$ are given by $X$ graded left $kX$ modules.
We denote the $kX$ action by $\tilde \la$, and  the grading by $\|.\|$. In a
representation of $D(X)$, the grading and $X$ action are related by
\[    \|y \tilde \la v\|\ =\  y\| v\| y^{-1} \quad y\in X\ ,v\in V\ ,
\]
and the action of $\delta_x\tens y\in D(X)$ is given by
\[
(\delta_x\tens y) \tilde \la v\ =\ \delta_{x,\|y \tilde \la v\|}\  y \tilde \la
v\ .
\]

We give an operation $\chi$ which sends representations of $D(H)$ to
representations of $D(X)$ in the following manner: Let $W$ be a
representation of $D(H)$, as described in the previous section.
Then as vector spaces, $\chi W$ is the same as $W$. The
$X$-grading $\|.\|$ on $\chi W$
is defined by
 $\|\chi(w)\|=\<w\>^{-1}|w|$, and the action of $us$ in $kX$
is given by
\[
us \tilde\la \chi(w)\ =\ \chi\bigg( \big((s\ra|w|^{-1})\la w\big)\ra u^{-1}
\bigg)
\ ,\quad s\in M,\ u\in G\ .
\]

\begin{propos} The map $\chi$ gives a 1-1 correspondence between
representations of $D(H)$ and $D(X)$.
\end{propos}
\proof
We prove that $\|.\|$ and $\tilde\la$ give a representation of $D(X)$.
First we must show that $\tilde\la$ is an action, which means that
\[vt \tilde\la \big(us \tilde\la \chi(w)\big)\ =\ vtus \tilde\la \chi(w)\ ,\]
for all $s,t\in M$ and $u,v\in G$. Note that
\begin{eqnarray*}
vt \tilde\la \big(us \tilde\la \chi(w)\big) & = & vt \tilde \la
\chi\bigg( \big((s\ra|w|^{-1})\la w\big)\ra u^{-1} \bigg) \\
& = &
\chi\bigg( \big((t\ra|\bar w|^{-1})\la \bar w\big)\ra v^{-1} \bigg)\ ,
\end{eqnarray*}
where $\bar w=\big((s\ra|w|^{-1})\la w\big)\ra u^{-1}$. On the other hand
 we have $vtus=v(t\la u)(t\ra u)s$, where $v(t\la u)\in G$ and
$(t\ra u)s\in M$, so
\begin{eqnarray*}
vtus \tilde\la \chi(w)& = &
\chi\bigg(\big(((t\ra u)s\ra|w|^{-1})\la w\big)\ra (t\la u)^{-1}v^{-1} \bigg)
\\
& = &\chi\bigg(\bigg(\big(((t\ra u)s\ra|w|^{-1})\la w\big)\ra (t\la u)^{-1}
\bigg)\ra v^{-1} \bigg)\ .
\end{eqnarray*}
We must therefore show that
\[
(t\ra|\bar w|^{-1})\la \bar w \ =\
\big(((t\ra u)s\ra|w|^{-1})\la w\big)\ra (t\la u)^{-1}\ .
\]
To do this we have to find the $G$ grade of $\bar w$, which we do as follows:
Put $\bar w=w'\ra u^{-1}$, where $w'=(s\ra|w|^{-1})\la w$.  Then using the
properties listed in 4.1,
\begin{eqnarray*}
|w'| & = & \big| (s\ra|w|^{-1})\la w \big| \ =\
(s\ra|w|^{-1})\la |w|\ =\ (s\la|w|^{-1})^{-1}\ , \\
\<w'\> & = & (s\ra|w|^{-1})\<w\> s^{-1}\ ,\\
|\bar w| & = & | w'\ra u^{-1}| \ =\ \big(\<w'\> \la u^{-1} \big)^{-1} |w'|
u^{-1}
 \ =\ \big(\<w'\> \la u^{-1} \big)^{-1} (s\la|w|^{-1})^{-1} u^{-1}\ .
\end{eqnarray*}
Putting $\bar t = t\ra u(s\la |w|^{-1})$,  and $\bar v=(t\la u)^{-1}$, we have
\begin{eqnarray*}
\bigg(\big( (t\ra u)s\ra |w|^{-1}\big) \la w\bigg)\ra (t\la u)^{-1} & = &
\bigg(
\big( (t\ra u)s\ra (s\la |w|^{-1})\big)  \la  \big( s\ra |w|^{-1}\big)
\la w \bigg)
\ra (t\la u)^{-1} \\
& = & \big( \bar t \la w'\big) \ra \bar v  \\
& = & \bigg( \bar t \ra \big(\<w'\>(\bar t \ra|w'|)^{-1}\la \bar v\big) \bigg)
\la  \bigg( w'\ra\big((\bar t \ra |w'|)^{-1}\la \bar v\big) \bigg) \ ,
\end{eqnarray*}
where the derivation of the last line uses 4.1(iv). To simplify the last line,
we calculate
\begin{eqnarray*}
\bar t \ra|w'| & = & \big(t\ra u(s\la|w|^{-1})\big)
\ra (s\la|w|^{-1})^{-1}\ =\ t\ra u\ ,\\
(\bar t \ra|w'|)^{-1} \la \bar v & = & (t\ra u)^{-1} \la (t\la u)^{-1}\ =\
u^{-1}\ .
\end{eqnarray*}
Then we have
\begin{eqnarray*}
\bigg(\big( (t\ra u)s\ra |w|^{-1}\big) \la w\bigg)\ra (t\la u)^{-1} & = &
\bigg( \bar t \ra \big(\<w'\> \la \bar u^{-1}\big) \bigg)
\la  \bigg( w'\ra u^{-1} \bigg) \\
& = & \big(\bar t\ra(\<w'\> \la \bar u^{-1})\big) \la \bar w \\
& = & \bigg(\big(t\ra u(s\la |w|^{-1})\big)
\ra\big(\<w'\>\la u^{-1}\big)\bigg)
 \la \bar w \\
& = & (t\ra|\bar w|^{-1})\la \bar w\ ,
\end{eqnarray*}
as required.
Next we must show that $\|y\tilde\la \chi(w)\|=y\|\chi(w)\|y^{-1}$.
This proceeds by a similar calculation, and is left to the reader!

Finally we show that $\chi$ is a 1-1 correspondence by giving its inverse
$\chi^{-1}$ in the following form.
Let $V$ be a representation of $D(X)$,
with $kX$ action $\tilde\la$, and $X$-grading $\|.\|$. Define a $D(H)$
representation as follows: As a vector space $\chi^{-1}V$ will be the
same as $V$. There are $G$ and $M$ gradings given by the factorisation
\[
\|v\| \ =\ \<\chi^{-1}v\>^{-1}|\chi^{-1}v|\ , \quad \<\chi^{-1}v\>\in M\ ,\quad
|\chi^{-1}v|\in G\ .
\]
The actions of $s\in M$ and $u\in G$ are given by
\[
s\la \chi^{-1}v\ =\ \chi^{-1}\big( (s\ra|\chi^{-1}v|) \tilde\la v \big)
\quad,\quad
\chi^{-1}v\ra u\ =\ \chi^{-1}\big( u^{-1} \tilde\la v \big) \ .
\]
\endproof

A 1-1 correspondence between the representations of $D(X)$ and $D(H)$
suggests an algebra isomorphism between $D(X)$ and $D(H)$.
We show that this is in fact the case.

\begin{propos}
There is an algebra isomorphism $\psi:D(H)\to D(X)$ defined by the formula
\[
\psi\big( \delta_s\tens u\tens t\tens \delta_v\big) \ =\
\delta_{u^{-1}s^{-1}(t\la v)u}\tens u^{-1}(t\ra v)\ ,
\]
and such that $\chi$ is induced by pull back along $\psi$.
\end{propos}
\proof (a)
To prove that $\psi$ is an algebra isomomorphism,
 we take $\alpha=\delta_x\tens q\tens t\tens\delta_v$, and
$\beta=\delta_s\tens u\tens r\tens\delta_w$, and their images
$\psi(\alpha)=\delta_{\bar s}\tens \bar u$ and
$\psi(\beta)=\delta_{\bar t}\tens \bar v$. Then the product in $D(X)$ is
\[
\psi(\alpha) \psi(\beta)\ =\ \delta_{\bar u^{-1} \bar s \bar u,\bar t}
\big( \delta_{\bar s}\tens \bar u\bar v \big)\ .
\]
Expanding this using the formula for $\psi$ gives
\[
\psi(\alpha) \psi(\beta)\ =\
\delta_{(t\ra v)^{-1}x^{-1}tv,u^{-1}s^{-1}(r\la w)u}
\big( \delta_{q^{-1}x^{-1}(t\la v)q}\tens q^{-1}(t\ra v)u^{-1}(r\ra w)\big)\ .
\]
Now using Proposition 3.2 to calculate $\alpha\beta$ in $D(H)$,
\[
\alpha \beta\ =\
\delta_{s'\ra u',x}\ \delta_{v',r\la w}\ \big(
\delta_{s'} \tens u'q \tens t'r \tens \delta_{w} \big)\ ,
\]
where
\[
t'=t\ra(s\la u)^{-1}\ ,\ v'=(s\la u) vu^{-1}\ ,\ s'=t's(t\ra vu^{-1})^{-1}\ ,\
u'=(t\ra vu^{-1})\la u\ .
\]
Applying $\psi$ to this product gives
\[
\psi(\alpha \beta)\ =\ \delta_{s'\ra u',x}\ \delta_{v',r\la w}\  \big(
\delta_{ (s'u'q)^{-1}(t'r\la w) u'q} \tens (u'q)^{-1}(t'r\ra w) \big)\ .
\]
First we establish the equality of the delta functions
\[
\delta_{(t\ra v)^{-1}x^{-1}tv,u^{-1}s^{-1}(r\la w)u}\ =\
\delta_{s'\ra u',x}\ \delta_{v',r\la w}\ .
\]

Firstly, suppose that the left hand side is non-zero. Then the uniqueness of
the factorisation $X=MG$ gives $(t\ra v)^{-1}x^{-1}t=(s\ra u)^{-1}$ and
$v=(s\la u)^{-1}(r\la w)u$. Rearranging these yields
$ r\la w=(s\la u)vu^{-1} $ and
$ x=t(s\ra u)(t\ra v)^{-1} $,
so the right hand side is also non-zero. Similarly for the other way around.

Next, assuming the above delta functions are not zero, we show that
\[
q^{-1}x^{-1}(t\la v)q
\ =\ (s'u'q)^{-1}(t'r\la w) u'q\ .
\]
This reduces to showing that
\[
s'u'x^{-1}(t\la v)
\ =\ (t'r\la w) u'\ .
\]
Since $r\la w=(s\la u)vu^{-1} $, we can calculate
\[
(t'r\la w) u'\ =\ \big(t\la (s\la u)^{-1}\big)^{-1}(t\la v)\ ,
\]
and since $ x=t(s\ra u)(t\ra v)^{-1} $,
\[
s'u'x^{-1}(t\la v)\ =\ \big(t\la (s\la u)^{-1}\big)^{-1}(t\la v)\ ,
\]
so we have the required identity.

Lastly, again under the assumption that the
delta functions are not zero, we show that
\[
q^{-1}(t\ra v)u^{-1}(r\ra w)
\ =\ (u'q)^{-1}(t'r\ra w)\ .
\]
This reduces to showing that
\[
(t\ra v)u^{-1}
\ =\ u'^{-1}(t'r\ra w)(r\ra w)^{-1}\ .
\]
On substituting for $r\la w$,
\[
u'^{-1}(t'r\ra w)(r\ra w)^{-1}\ =\
u'^{-1}\big(t'\ra (r\la w)\big) (r\ra w)(r\ra w)^{-1}\ =\
u'^{-1}(t\ra vu^{-1})\ =\ (t\ra v)u^{-1}\ .
\]

This completes the proof that $\psi$ is an algebra map.
To show that it is an isomorphism, we can give an
explicit form for the inverse. For $s,t\in M$ and $u,v\in G$, we have
\[
\psi^{-1}\big( \delta_{su}\tens tv \big)\ =\
\delta_{(s^{-1}\ra(t\la v))^{-1}}\tens (t\la v)^{-1} \tens
(t\ra v\alpha^{-1}v)\tens \delta_{v^{-1}\alpha}\ ,
\]
where $\alpha=t^{-1}\la u^{-1}(s^{-1}t\la v)$.
\endproof
\proof (b): We next show that the maps $\psi$ and $\chi$ are linked
 by the equation
 $\chi\big((a\tens h)\la w\big)=\psi(a\tens h)\tilde\la \chi(w)$.

Begin by
\begin{eqnarray*}
\chi\big((\delta_s\tens u \tens t\tens\delta_v)\la w\big)& =&
\delta_{s,\<t\la w\>} \delta_{v,|w|}
\chi\big((t\la w)\ra u\big)\ ,\\
\psi(\delta_s\tens u \tens t\tens\delta_v)\tilde \la \chi w & = &
\big( \delta_{u^{-1}s^{-1}(t\la v)u}\tens u^{-1}(t\ra v) \big)
\tilde \la \chi w \\
& = &\delta_{u^{-1}s^{-1}(t\la v)u,\|u^{-1}(t\ra v)\tilde \la \chi w\|}.
u^{-1}(t\ra v)\tilde \la \chi w \\
& = &
\delta_{u^{-1}s^{-1}(t\la v)u,u^{-1}(t\ra v) \|\chi w\|(t\ra v)^{-1}u}.
 \chi \bigg( \big( ( (t\ra v)\ra|w|^{-1})\la w \big)\ra u \bigg)
\end{eqnarray*}
Now we simplify the delta function: If the delta function is non-zero, then
\[
s^{-1}t v\ =\
s^{-1}(t\la v)(t\ra v)\ =\ (t\ra v) \|\chi w\|\ =\ (t\ra v) \< w\>^{-1}|w|\ .
\]
By the uniqueness of factorisation, we find that
$s^{-1}t=(t\ra v) \< w\>^{-1}$ and $v=|w|$. But then $s=t\< w\>
(t\ra |w|)^{-1}=\<t\la w\>$, giving the equality we want.
\endproof

We have now discussed the algebra structure and algebra representations.
To examine the coalgebra structure we turn to tensor
 products of representations.
Recall that for a general Hopf algebra $H$ with representations $V$
 and $V'$, the tensor product representation $V\tens V'$ is defined by
the action
\[
h\la (v\tens v')\ =\ \sum h_{(1)}\la v \tens h_{(2)}\la v'\ .
\]
In the case of representations
$V$ and $V'$ of $D(X)$, this formula leads to the equations
\[
  y\tilde\la (v\tens v')\ =\ y\tilde\la v \tens y
\tilde\la v' \quad\text{and}\quad
  \| v \tens  v'\|\ =\ \| v \| \| v'\|\ .
\]
If we take representations $W$ and $W'$ of $D(H)$, we get
the following actions and gradings on $W\tens W'$:
\begin{eqnarray*}
t\la (w\tens w') & = & (t\la w) \tens \big( (t\ra |w|)\la w'\big)  \quad
t\in M\ , \\
|w\tens w'| & = & |w||w'| \ ,\\
(w\tens w') \ra u & = & \big( w\ra (\<w'\>\la u)\big) \tens (w'\ra u)  \quad
u\in G\ , \\
\<w\tens w'\> & = & \<w\>\<w'\>   \ .
\end{eqnarray*}
We find that $\chi$ does not preserve tensor products of representations.
To correct for this we introduce another map
$c:\chi W \tens \chi W'\to \chi(W\tens W')$, defined by
\[
c\big( \chi(w)\tens\chi(w')\big)\ =\
\chi\big((\<w'\>\ra |w|^{-1})\la w\tens w'\big)\ .
\]

\begin{propos}
The map
$c$ is a $D(X)$ module map, i.e.
it preserves the grading and action in the following manner:
\[
  \| c\big( \chi w\tens\chi w' \big) \|\ =\  \|  \chi w\tens\chi w'  \|
\]
\[
y \tilde\la  c\big( \chi w\tens\chi w' \big)\ =\ c\big( y \tilde\la
\big( \chi w\tens\chi w' \big)\big)\ .
\]
for all $y\in X, w\in W$ and $w'\in W'$.
\end{propos}
\proof
First we begin with the grading. We know that
\[
 \|  \chi w\tens\chi w'  \|\ =\  \|  \chi w  \|  \|  \chi w'  \|\ =\
\<w \>^{-1} |w|  \<w' \>^{-1} |w'|\ .
\]
On the other hand
\begin{eqnarray*}
  \| c\big( \chi w\tens\chi w' \big) \| & = &
  \| \chi\big((\<w'\>\ra |w|^{-1})\la w\tens w'\big) \| \\
& = & \<(\<w'\>\ra |w|^{-1})\la w\tens w' \>^{-1}
 \big|(\<w'\>\ra |w|^{-1})\la w\tens w'\big|  \\
 & = &  \<w' \>^{-1}   \<(\<w'\>\ra |w|^{-1})\la w \>^{-1}
 \big|(\<w'\>\ra |w|^{-1})\la w\big|   \big|w'\big| \ .
\end{eqnarray*}
By 4.1 we can calculate
\begin{eqnarray*}
\<(\<w'\>\ra |w|^{-1})\la w \> & = &
\big(\<w'\>\ra |w|^{-1}\big) \<w\> \big((\<w'\>\ra |w|^{-1})\ra |w|\big)^{-1}
\\
& = &
\big(\<w'\>\ra |w|^{-1}\big) \<w\> \<w'\>^{-1} \\
\big|(\<w'\>\ra |w|^{-1})\la w \big| & = &  \big(\<w'\>\ra |w|^{-1}\big) \la
|w|
\ =\ \big(\<w'\>\la |w|^{-1}\big)^{-1}\ ,
\end{eqnarray*}
which gives the result.

Now we shall show that $c$ preserves the action,
which we do seperately for the two factors $G$ and $M$ of $X$.

For the $G$ action, we must show that
\[
c\big(u\tilde\la (\chi w\tens\chi w')\big)\ =\
u\tilde\la  c\big(\chi w\tens\chi w'\big)\ .
\]
By the definitions
\begin{eqnarray*}
u\tilde\la (\chi w\tens\chi w') & = &
\chi (w\ra u^{-1})\tens\chi (w'\ra u^{-1})\ ,\\
c\big(u\tilde\la (\chi w\tens\chi w')\big) & = &
\chi\bigg(
\big( (\<w'\ra u^{-1}\>\ra| w\ra u^{-1} |^{-1})\la (w\ra u^{-1})\big)
\tens(w'\ra u^{-1})\bigg)\ .
\end{eqnarray*}
Using the properties of the $G$ and $M$ gradings,
\begin{eqnarray*}
\<w'\ra u^{-1}\>\ra| w\ra u^{-1} |^{-1} & = &
\big( \<w'\> \ra |w|^{-1} \big)\ra \big( \<w\> \la u^{-1} \big) \ ,\\
\big( \<w'\ra u^{-1}\>\ra| w\ra u^{-1} |^{-1} \big)\la (w\ra u^{-1}) & = &
\big( \big( \<w'\> \ra |w|^{-1} \big)\ra \big( \<w\> \la u^{-1} \big)
 \big)\la (w\ra u^{-1}) \\
& = &
\big( (\<w'\> \ra |w|^{-1}) \la w \big)
\ra\big(\<w'\> \la u^{-1}\big)\ ,
\end{eqnarray*}
so that
\[
c\big(u\tilde\la (\chi w\tens\chi w')\big)\ =\
\chi\bigg(
\big( ( (\<w'\> \ra |w|^{-1})\la w) \ra (\<w'\>\la u^{-1}) \big)
\tens(w'\ra u^{-1})
\bigg)\ .
\]
Next we calculate
\begin{eqnarray*}
u\tilde\la  c\big(\chi w\tens\chi w'\big) & = &
u\tilde\la  \chi\bigg( \big( (\<w'\> \ra |w|^{-1}) \la w\big)
\tens w'\bigg) \\
& = & \chi\bigg(
\big( \big( (\<w'\> \ra |w|^{-1}) \la w\big)
\tens w'\big) \ra u^{-1} \bigg) \\
& = & \chi\bigg(
\big( ( (\<w'\> \ra |w|^{-1})\la w) \ra (\<w'\>\la u^{-1}) \big)
\tens(w'\ra u^{-1})
\bigg)\ ,
\end{eqnarray*}
which agrees with the alternative calculation.

Lastly we shall show that $c$ preserves the $M$ action. We have
\[
s\tilde\la (\chi w\tens\chi w') \ = \
\chi \big((s\ra |w|^{-1})\la w\big)\tens
\chi \big((s\ra |w'|^{-1})\la w'\big)   \ .
\]
Using the equations
\begin{eqnarray*}
\<(s\ra |w'|^{-1})\la w'\> & = & \big(s\ra |w'|^{-1}\big)\<w'\>s^{-1}  \\
\big|(s\ra |w|^{-1})\la w\big| & = & \big(s\ra |w|^{-1}\big)\la |w|\ =\
(s\la |w|^{-1})^{-1}\ ,
\end{eqnarray*}
we see that
\begin{eqnarray*}
c\big( s\tilde\la (\chi w\tens\chi w') \big) & = &
\chi\bigg(
\big(
\big(
\< (s\ra |w'|^{-1})\la w' \>
\ra \big| (s\ra |w|^{-1})\la w \big|^{-1}
\big)
\la\big( (s\ra |w|^{-1})\la w\big)
\big)
\tens \big( (s\ra |w'|^{-1})\la w' \big)
\bigg)   \\
 & = &
\chi\bigg(
\big(
\big(
(s\ra |w'|^{-1})\<w'\>s^{-1}
\ra
(s\la |w|^{-1})
\big)
 (s\ra |w|^{-1})\la w
\big)
\tens \big( (s\ra |w'|^{-1})\la w' \big)
\bigg)   \\
 & = &
\chi\bigg(
\big(
\big(
(s\ra |w'|^{-1})\<w'\>\ra|w|^{-1}
\big)
\la w
\big)
\tens \big( (s\ra |w'|^{-1})\la w' \big)
\bigg)  \ .
\end{eqnarray*}
On the other hand,
\begin{eqnarray*}
s\tilde \la c\big(\chi w\tens\chi w'\big) & = &
s\tilde \la \chi \big(\bar w \tens w'\big)  \\
 & = & \chi \big( (s\ra|\bar w \tens w'|^{-1})\la (\bar w \tens w') \big) \\
& = & \chi \big( (\bar t\la \bar w) \tens ((\bar t\ra|\bar w|)\la w') \big) \ ,
\end{eqnarray*}
where $\bar w=\big(\<w'\>\ra| w |^{-1}\big)\la w$ and $\bar t =
s\ra|\bar w \tens w'|^{-1}$. Next we calculate
\begin{eqnarray*}
|\bar w \tens w'| & = & |\bar w | | w'|\ ,\\
|\bar w | & = & \big( \<w'\>\la| w |^{-1} \big)^{-1}\ ,\\
\bar t \la \bar w  & = &
\big( (s\ra  | w' |^{-1} )\<w'\> \ra | w |^{-1} \big) \la w\ ,\\
(\bar t \ra |\bar w|) \la w' & = & (s \ra | w'|^{-1}) \la w'\ ,
\end{eqnarray*}
and substitution of these gives the answer.
\endproof

We find that
$(\chi,c)$ is a monoidal functor between the categories
of representations.
By Tannaka-Krein reconstruction theory, we then expect that the Hopf
algebras
$D(H)$ and $D(X)$
(which are already isomorphic as algebras) are related by a cocycle $F$
in the following manner:

If $F\in D(X)\tens D(X)$ is an invertible, we define a new coproduct
$\tilde \Delta$  for $D(X)$ and element $R$ by the formula \cite{Dri:2}
\[
\tilde \Delta(h)\ =\ F(\Delta h)F^{-1}\ ,\quad\text{and}\quad
\tilde R\ =\ (\tau F)RF^{-1}\ .
\]
This also gives a quasitriangular Hopf algebra if $F$ is a 2-cocycle. The
condition for a 2-cocycle is
$(1\tens F)(id\tens \Delta)F=(F\tens 1)(\Delta\tens id)F$.
The new Hopf algebra has the same representations as $D(X)$ with a different
tensor product, i.e.\ is monoidally equivalent to it.

\begin{propos} The following  formula defines a 2-cocycle:
\[
F\ =\ \sum_{x\in X\ ,t\in M\ ,v\in G}
 \delta_x\tens t^{-1}\tens \delta_{tv}\tens e\ .
\]
\end{propos}
\proof
We follow our usual convention that $r$, $s$ and $t$ are in $M$, that
$u$, $v$ and $w$ are in $G$, and that $x$, $y$, $a$, $b$ and $z$ are in $X$.
First we calculate
\begin{eqnarray*}
(id\tens \Delta)F & = & \sum_{x , t , v}
\sum_{a,b:ab=tv}
 \delta_x\tens t^{-1}\tens \delta_{a}\tens e\tens \delta_{b}\tens e\ ,  \\
(1\tens F)(id\tens \Delta)F & = & \sum_{x , t , v}
\sum_{a,b:ab=tv}\sum_{y , s , u}
\delta_x\tens t^{-1}\tens
\big( \delta_y\tens s^{-1} \big)\big( \delta_a\tens e \big) \tens
\big( \delta_{su}\tens e \big)\big( \delta_b\tens e \big)  \\
& = & \sum_{x , t , v}
\sum_{a,b:ab=tv}\sum_{y , s , u}
\delta_{sys^{-1},a}\delta_{su,b}.
\delta_x\tens t^{-1}\tens \delta_y\tens s^{-1} \tens \delta_{su}\tens e \ ,\\
& = & \sum_{x , t , v, s , u}
\delta_x\tens t^{-1}\tens \delta_{s^{-1}tvu^{-1}}\tens s^{-1} \tens
\delta_{su}\tens e \ .
\end{eqnarray*}
On the other hand,
\begin{eqnarray*}
(\Delta\tens id)F & = & \sum_{\bar x ,\bar  t ,\bar  v}
\sum_{a,b:ab=\bar x }
 \delta_a\tens\bar  t^{-1}
\tens \delta_{b}\tens \bar  t^{-1}\tens
 \delta_{\bar  t \bar  v}\tens e\ ,  \\
(F\tens 1)(\Delta\tens id)F & = & \sum_{\bar x ,\bar  t ,\bar  v,z,r,w}
\sum_{a,b:ab=\bar x }
\big( \delta_{z}\tens r^{-1} \big)  \big( \delta_a\tens\bar  t^{-1} \big)
\tens
\big(  \delta_{rw}\tens e\big)  \big( \delta_{b}\tens \bar  t^{-1} \big)
 \tens
 \delta_{\bar  t \bar  v}\tens e\ ,  \\
& = & \sum_{\bar x ,\bar  t ,\bar  v,z,r,w}
\sum_{a,b:ab=\bar x }
\delta_{rzr^{-1},a} \delta_{rw,b}.
\delta_{z}\tens r^{-1} \bar  t^{-1} \tens  \delta_{rw}\tens \bar  t^{-1}
\tens  \delta_{\bar  t \bar  v}\tens e\ ,  \\
& = & \sum_{\bar x ,\bar  t ,\bar  v,r,w}
\delta_{r^{-1}\bar x w^{-1}}\tens r^{-1} \bar  t^{-1} \tens  \delta_{rw}\tens
\bar  t^{-1}  \tens  \delta_{\bar  t \bar  v}\tens e\ .
\end{eqnarray*}
Now if we put $\bar t=s$, $\bar v=u$, $r=s^{-1}t$, $w=vu^{-1}$, and
$\bar x=s^{-1}txvu^{-1}$, we see that the two expressions are equal.
\endproof

We now know that twisting the coproduct of $D(X)$ by $F$ gives another
quasi-triangular Hopf algebra. It remains to relate this to the $D(H)$
Hopf algebra in the following manner:

\begin{propos}
\[
\tilde\Delta(h)\ =\ F.(\Delta_{D(X)} h).F^{-1}\ =\
(\psi\tens\psi)\Delta_{D(H)}(\psi^{-1}h)\
\]
for all $h\in D(X)$. Hence the the twisting of $D(X)$ by $F$ is isomorphic to
$D(H)$.
\end{propos}
\proof
Here we use the conventions that $x,y,a,b\in X$, $r,s,t,\bar t,a',b'\in M$, and
$u,v,w,\bar v,x',y'\in G$.
Begin by calculating
\[
\delta_x\tens y\ =\ \psi(\delta_s\tens u\tens t\tens\delta_v)\ =\
\delta_{u^{-1}s^{-1}(t\la v)u}\tens u^{-1}(t\ra v)\ .
\]
Then we have
\begin{eqnarray*}
(\psi\tens\psi)\Delta(\delta_s\tens u\tens t\tens\delta_v)& =&
\sum_{a'b'=s,x'y'=v}
\psi\big(\delta_{a'}\tens (b'\la u)\tens t\tens\delta_{x'}\big)
\tens
\psi\big(\delta_{b'}\tens u\tens (t\ra x')\tens\delta_{y'}\big) \\
 =\   \sum_{a'b'=s,x'y'=v} & &
\delta_{ (b'\la u)^{-1}a'^{-1}(t\la x')(b'\la u) }
\tens  (b'\la u)^{-1} (t\ra x') \tens
\delta_{ u^{-1}b'^{-1}\big(((t\ra x'))\la y'\big)u }
\tens u^{-1}(t\ra v) \\
 =\   \sum_{a'b'=s,x'y'=v} & &
\delta_{ (b'\la u)^{-1}a'^{-1}(t\la x')(b'\la u) }
\tens  (b'\la u)^{-1} (t\ra x') \tens
\delta_{ \bar t \bar v }
\tens y
\end{eqnarray*}
where we have used the relabelings
\[
\bar t\ =\ (b'\ra u)^{-1}\ ,\ \bar v\ =\ (b'\la u)^{-1}
(t\la x')^{-1}(t\la v)u\ .
\]
Then we can find
\begin{eqnarray*}
b'\la u & = & (\bar t^{-1} \la u^{-1} )^{-1}\ ,\\
t\la x'  & =& (t\la v)u\bar v^{-1}(b'\la u)^{-1}\ =\
(t\la v)u\bar v^{-1}(\bar t^{-1} \la u^{-1} )\ .
\end{eqnarray*}
Now use the relation $a'b'=s$ to show that
\begin{eqnarray*}
(b'\la u)^{-1}a'^{-1}(t\la x')(b'\la u) & = &
(\bar t^{-1} \la u^{-1} )a'^{-1}(t\la v)u\bar v^{-1} \\
 & = &
(\bar t^{-1} \la u^{-1} )(\bar t^{-1} \ra u^{-1} )s^{-1}(t\la v)u\bar v^{-1} \\
 & = &
\bar t^{-1}  u^{-1}s^{-1}(t\la v)u\bar v^{-1} \ =\ \bar t^{-1} x \bar v^{-1}\ .
\end{eqnarray*}
Also we consider
\begin{eqnarray*}
(b'\la u)^{-1} (t\ra x') & = & (\bar t^{-1} \la u^{-1} )  (t\la x')^{-1} tx' \\
 & = &  \bar v u^{-1} (t\la v)^{-1} tx' \\
 & = &  \bar v u^{-1} (t\ra v)v^{-1} x'\ =\  \bar v yv^{-1} x'\ .
\end{eqnarray*}
On the other hand, we calculate
\begin{eqnarray*}
F . \Delta (\delta_x\tens y) & = &
\bigg(\sum_{\bar x,\bar t,\bar v}
 \delta_{\bar x}\tens \bar t^{-1}\tens \delta_{\bar t\bar v}\tens  e \bigg)
\bigg(\sum_{ab=x}
 \delta_{a}\tens y\tens \delta_{b}\tens  y \bigg)  \\
& = & \sum_{\bar t,\bar v}
\delta_{\bar t^{-1} x \bar v^{-1} }\tens \bar t^{-1}y\tens
\delta_{\bar t\bar v}\tens y\ .
\end{eqnarray*}
The next stage is to calculate
\begin{eqnarray*}
F . \Delta (\delta_x\tens y).F^{-1} & = &
\bigg( \sum_{\bar t,\bar v}
\delta_{\bar t^{-1} x \bar v^{-1} }\tens \bar t^{-1}y\tens
\delta_{\bar t\bar v}\tens y \bigg)
\bigg( \sum_{z,r,w} \delta_{z}\tens r\tens \delta_{rw}\tens e \bigg)  \\
& = & \sum_{r,w,\bar t,\bar v:rw=y^{-1}\bar t\bar v y}
\delta_{\bar t^{-1} x \bar v^{-1}}\tens \bar t^{-1} yr
 \tens \delta_{\bar t\bar v}\tens y\ .
\end{eqnarray*}
Then the equation $rw=y^{-1}\bar t\bar v y$ can be rearranged to give
$\bar t^{-1} yr=\bar v y w^{-1}$, which gives the required equality of the two
expressions (with $w=x'^{-1}v$).
\endproof

This is also true at the level of quasitriangular structures. Recall
first that if $R$ is a quasitriangular structure, then so is $\tau R^{-1}$.
In our above conventions we have:

\begin{propos}
The relation between the quasi-triangular structures on $D(X)$ and $D(H)$
is given by the formula
$(\tau F)(\tau R^{-1}) F^{-1}=\psi\tens\psi \big(R_{D(H)} \big)$.
\end{propos}
\proof
First we calculate
\begin{eqnarray*}
(\tau F)R & = &
\bigg(\sum_{x,\bar t,\bar v} \delta_{\bar t\bar v}\tens  e
\tens  \delta_{x}\tens \bar t^{-1} \bigg)
\bigg(\sum_{y,z} \delta_{y}\tens e \tens  \delta_{z}\tens y \bigg) \\
 & = &
\sum_{x,\bar t,\bar v} \delta_{\bar t\bar v}\tens  e
\tens  \delta_{x}\tens \bar v \ ,\\
(\tau F)R F^{-1} & = &
\bigg(
\sum_{x,\bar t,\bar v} \delta_{\bar t\bar v}\tens  e
\tens  \delta_{x}\tens \bar v
 \bigg)
\bigg(\sum_{y,r,w} \delta_{y}\tens  r
\tens  \delta_{rw}\tens e \bigg) \\
& = &
\sum_{x,\bar t,\bar v} \delta_{\bar t\bar v}\tens  r
\tens  \delta_{x}\tens \bar v\ ,
\end{eqnarray*}
where $r\in M$ is the solution to the factorisation $rw=\bar v^{-1} x \bar v$.

On the other hand,
\begin{eqnarray*}
(\psi\tens\psi)R_{D(H)}
& = & (\psi\tens\psi)\sum_{s'\in M, \ u'\in G}
\delta_{s'}\tens u'\tens e\tens 1\tens 1\tens e \tens s'\tens \delta_{u'} \\
& = & \sum_{s',t'\in M, \ u',v'\in G}
\delta_{u'^{-1}s'^{-1}v'u'}\tens  u'^{-1}
\tens  \delta_{t'^{-1}(s'\la u')}\tens (s'\ra u')\ .
\end{eqnarray*}
We can reorganise this by using $s=(s'\ra u')^{-1}$, $u=(s'\la u')^{-1}$,
$v=uv'u'$ and $t=t'^{-1}$ to get
\[
(\psi\tens\psi)R_{D(H)}\ =\ \sum_{s,t\in M, \ u,v\in G}
\delta_{sv}\tens  (s\la u)
\tens  \delta_{tu^{-1}}\tens s^{-1}\ .
\]
Now if we calculate (left to the reader), we find
\[
\big(\tau (\psi\tens\psi)R_{D(H)}\big)(\tau F)R F^{-1}\ =\ 1_{D(X)}\tens
1_{D(X)}\ .
\]
\endproof

The cocycle $F$ is not a coboundary. This means that there is no inner
automorphism relating the two coproduct structures on $D(X)$, although the
Hopf algebras may still be isomorphic via an outer
automorphism. To see this, we note
that if $F$ was a coboundary, there would be an invertible element
$\gamma\in D(X)$ so that $\gamma\tens\gamma = F.\Delta \gamma$ in $D(X)$.
We shall take $\gamma$ of the form
\[
\gamma\ =\ \sum_{x,y\in X}\gamma_{x,y} \delta_x\tens y\ ,
\]
and $F$ of the form given previously. Then we can work out the following
expressions as
\begin{eqnarray*}
F.\Delta \gamma& = & \sum_{z,t,v,y} \gamma_{tzv,y}
\delta_z \tens t^{-1}y \tens  \delta_{tv} \tens y\ , \\ \gamma \tens \gamma
& = & \sum_{z,t,v,y,x} \gamma_{z,x} \gamma_{tv,y}
\delta_z \tens x \tens  \delta_{tv} \tens y\ .
\end{eqnarray*}
If these are equal, we must have
\[
\gamma_{z,x} \gamma_{tv,y} \ =\
\left\{ \begin{array}{ll}
0
& x\neq t^{-1}y   \\
\gamma_{tzv,y}
&  x= t^{-1}y
\end{array} \right.
\]
To make these equations clearer, we put
\[
\beta_{t^{-1}y,v,y}\ =\ \gamma_{tv,y} \ ,
\]
and write $z=su$ to rewrite the previous equation as
\[
\beta_{s^{-1}x,u,x} \beta_{t^{-1}y,v,y} \ =\
\left\{ \begin{array}{ll}
0
& x\neq t^{-1}y   \\
\beta_{s^{-1}x,uv,y}
&  x= t^{-1}y
\end{array} \right.\ .
\]
Now suppose that there are $x_0$, $s_0$, and
$u_0$ with $\beta_{s_0^{-1}x_0,u_0,x_0}\neq 0$. The equation then shows
that
$\beta_{t^{-1}y,v,y}=0$ if $y\neq s_0^{-1}x_0$. If we consider the
expression $\beta_{s_0^{-1}x_0,u_0,x_0}\beta_{s_0^{-1}x_0,u_0,x_0}\neq 0$
in the equation, we find that $s_0=e$. The only possible non-zero
$\beta$s are of the form $\beta_{x_0,u,x_0}$, which corresponds to the
only possible non-zero
$\gamma$s being of the form $\gamma_{u,x_0}$. But this means that
$\gamma$
cannot be invertible
in $D(X)$ (unless one of the groups in the factorisation of $X$
is just the identity). We see that the non-abelian cohomology
$H^2(D(X))$ is non-trivial.


\end{document}